\documentstyle[epsfig,amsmath,amssymb,amsfonts,MATHYPH]{article}
\evensidemargin\oddsidemargin
\newcommand{\PSbox}[3]{\mbox{\rule{0in}{#3}\includegraphics{#1}\hspace{#2}}} 
\newcommand{\1}{{1\!\!1}}

\newcommand{\sevafig}[2]{\begin{figure}[h]\centerline{
 \epsfig{file=#1,width=#2}}
\bigskip\end{figure}}

\begin{document}

\newcommand{\gothk}{\mathfrak{k}}
\newcommand{\gothg}{\mathfrak{g}}

\large
                                                                    
\centerline{\Large\bf Sur l'homologie des espaces de n\oe uds non-compacts}

\vspace{2mm}

\centerline{\bf Victor Tourtchine}

\vspace{3mm}

%\hangindent=2cm

\begin{tabular}{p{1.3cm}p{12cm}}
&
{\small
{\bf Mots cl\'es.} Discriminant de l'espace de n\oe uds non-compacts,
big\`ebre des diagrammes de chordes, complexe de Hochschild,
 op\'erades des alg\`ebres de Poisson, 
de Gerstenhaber, de Batalin-Vilkovissky.}
\\
\end{tabular}

\vspace{2mm}

\begin{tabular}{p{1.3cm}p{12cm}}
&

{\small {\bf R\'esum\'e.} La suite
spectrale de Vassiliev, voir [V1], calcule l'homologie de l'espace de n\oe uds 
non-compacts dans ${\Bbb R}^d$, $d\ge 3$.
Dans ce travail on d\'ecrit le premier terme 
de cette suite spectrale
 en terme de l'homologie du 
complexe de 
Hochschild pour l'op\'erade des alg\`ebres de Poisson, si $d$ est impair 
(resp. pour l'op\'erade des alg\`ebres de Gerstenhaber, si $d$ est pair).
En particulier, la big\`ebre des diagrammes de cordes appara\^\i t comme
sous-espace de cette homologie. L'homologie en question peut \^etre vue comme 
l'espace des classes charact\'eristiques de l'homologie des complexes de 
Hochschild pour les alg\`ebres de Poisson (resp. de Gerstenhaber), \'etudi\'ees
comme alg\`ebres associatives. On donne aussi une simplification des 
calculs du premier terme de la suite spectrale de Vassiliev.}
\\
\end{tabular}

\vspace{4mm}

\noindent{\large\bf 0. Introduction}                         

\vspace{4mm}

\noindent{\bf 0.1. Histoire du sujet}

\vspace{2mm}

On appelle  {\it n\oe uds non-compacts} les applications lisses injectives
non-singuli\`eres ${\Bbb R}\hookrightarrow {\Bbb R}^d$, qui co\"\i ncident
hors d'un certain sous-ensemble compact de 
${\Bbb R}$ avec une application 
lin\'eaire fix\'ee. Les n\oe uds non-compacts forment un sous-ensemble ouvert 
partout dense dans l'espace (affine) ${\cal K}$ de toutes les applications 
lisses
${\Bbb R}\rightarrow {\Bbb R}^d$ avec le m\^eme comportement \`a 
l'infini. Son compl\'ement $\Sigma \subset {\cal K}$ est un {\it espace 
discriminant}, qui consiste des applications ayant des auto-intersections
ou des singularit\'es. Toute classe de cohomologie 
$\gamma \in H^i$ $({\cal K}\backslash \Sigma )$ de l'espace
de n\oe uds peut \^etre r\'ealis\'ee comme un indice d'enlacement
avec une cha\^\i ne convenable dans $\Sigma$ de codimension
$i+1$ dans ${\cal K}$.

Pour  simplifier (suivant [V5]), on va supposer que l'espace $\cal K$
est d'une dimension $\omega $
tr\`es grande mais finie.
L'explication rigoureuse de l'hypoth\`ese utilise des approximations
de dimension finie de l'espace ${\cal K}$ (voir [V1]).
Ci-dessous nous mettons entre guillemets les affirmations non-rigoureuses,
qui utilisent cette hypoth\`ese et qui sont \`a pr\'eciser.

L'outil principal de cette approche pour calculer l'homologie 
de l'espace de n\oe uds est la r\'esolution simpliciale (construite dans
[V1]) du discriminant $\Sigma$ appel\'ee discriminant 
r\'esolu $\sigma$. La projection naturelle 
$\pi : {\bar \sigma}\to {\bar \Sigma}$ induit 
une ``\'equivalence homotopique'' des compactifi\'es par un point
des espace $\sigma$ et $\Sigma$. Par la ``dualit\'e d'Alexander''  les
groupes d'homologie $\tilde{H}_*(\bar{\sigma})\equiv
\tilde{H}_*(\bar{\Sigma} )$  de ces compactifi\'es ``co\"\i ncident'' (\`a un 
d\'ecalage de dimension pr\`es) avec les groupes de cohomologie de
l'espace des n\oe uds:
$$
\tilde{H}^i({\cal K}\backslash \Sigma ,\Bbbk )\simeq\tilde{H}_{\omega -i-1}
(\bar{\Sigma},\Bbbk )\equiv\tilde{H}_{\omega -i-1}(\bar{\sigma}, \Bbbk ).
\eqno(0.1.1)$$
($\Bbbk$ est un anneau commutatif de coefficients.)

L'espace $\sigma$ poss\`ede une filtration naturelle

$$
\varnothing =\sigma_0\subset\sigma _1\subset\sigma_2\subset \dots .
\eqno(0.1.2)$$

{\bf Conjecture 0.1.3.} {\it La suite spectrale (de Vassiliev) 
associ\'ee \`a la filtration (0.1.2) et calculant l'homologie de Borel-Moore
 de la r\'esolution $\sigma$
d\'eg\'en\`ere sur $\Bbb Q$ au premier terme.} $\Box$

Il existe aussi une autre cojecture plus forte.

{\bf Conjecture 0.1.4.} {\it La filtration (0.1.2) se scinde homotopiquement,
autrement dit  $\bar{\sigma}$ est ``homotopiquement \'equivalent'' au 
bouquet  $
V_{i=1}^{+\infty}(\bar{\sigma}_i/\bar{\sigma}_{i-1}).$} $\Box$

Cette conjecture entra\^\i nerait 
la d\'eg\'en\'er\'escence de notre {\it suite spectrale principale} au
premier terme sur n'importe quel anneau commutatif $\Bbbk$ de
coefficients. Modulo cette conjecture, pour calculer la cohomologie de
l'espace de n\oe uds non-compacts dans ${\Bbb R}^d$, avec $d\ge 4$
(si $d=3$, la suite spectrale \'etudi\'ee ne calcule qu'un certain
sous-groupe dans la cohomologie de l'espace de n\oe uds), il suffit
de savoir le premier terme. 

D'un autre 
c\^ot\'e, dans les termes
$\sigma_i\backslash\sigma_{i-1}$  de la filtration on a une  
d\'ecomposition cellulaire tr\`es simple, qui ne d\'epend
que de la parit\'e de la dimension $d$ de l'espace ${\Bbb R}^d$
\`a un d\'ecalage de dimenson pr\`es). Ceci  rend 
trivial du point de vue g\'eom\'etrique le calcul du premier terme de la 
suite spectrale de Vassiliev.

Pour calculer le premier terme V.A.Vassiliev a introduit une {\it
filtration auxiliaire} sur les termes $\sigma_i\backslash
\sigma_{i-1}$. La suite spectrale associ\'ee \`a cette filtration
d\'eg\'en\`ere au deuxi\`eme terme, parce que son premier terme (pour
tout $i$) est concentr\'e en une seule ligne. Le z\'eroi\`eme terme de
la suite auxiliaire avec sa z\'eroi\`eme diff\'erentielle est une
somme directe de produits tensoriels de complexes des graphes
connexes. L'homologie du complexe des graphes connexes sur $m$
points est concentr\'ee dans une seule dimension et est isomorphe \`a
${\Bbb Z}^{(m-1)!}$, voir [V3], [V4].
  Elle a une description simple comme espace
engendr\'e par les arbres et quotient\'e par les relations de trois
termes, voir [V3], [T].  

Pour $d=3$ en cohomologie de degr\'e z\'ero la suite spectrale de Vassiliev 
calcule une partie de la cohomologie de degr\'e z\'ero ---
les {\it invariants} 
{\it de type fini}, que l'on 
peut d\'efinir de mani\`ere plus simple et g\'eom\'etrique, voir [ChDL].
L'objet dual \`a l'espace des invariants de type fini est la {\it big\`ebre des
diagrammes de cordes}, qui a \'et\'e intens\'ement \'etudi\'ee pendant les 
derni\`eres ann\'ees, voir [BN], [ChD], [ChDL], [K1], [Kn], [L], [NS], [S], 
[Z]. Le but de ce 
travail est de bien montrer que dans l'homologie sup\'erieure des
espaces de n\oe uds (non-compacts) on a aussi de tr\`es belles 
math\'ematiques.

\vspace{4mm}

\noindent{\bf 0.2. Contenu. R\'esultats principaux}

\vspace{2mm}

La cohomologie du complexe des graphes connexes
(sur $m$ points) ---
l'espace dual \`a celui que l'on consid\'erait dans la section pr\'ec\'edente,
a aussi une description tr\`es simple.

Consid\'erons une alg\`ebre de Lie sur ${\Bbb Z}$ libre avec $m$
g\'en\'erateurs. Consid\'erons son sous-espace lin\'eairement engendr\'e par 
les crochets tels, que chaque g\'en\'erateur y est pr\'esent\'e 
exactement une fois. Ce sous-espace est isomorphe \`a ${\Bbb Z}^{(m-1)!}$.
Il se trouve que la cohomologie en question est 
exactement ce sous-espace.

Cet isomorphisme vient de la construction suivante:

Consid\'erons l'espace des applications injectives d'un ensemble fini
$M$ de $m$ \'el\'ements dans ${\Bbb R}^d$, $d\ge 1$. Cet espace peut \^etre
vu comme un analogue de dimension finie de l'espace des n\oe uds.
Le discriminant correspondant (qui se compose des applications non-injectives)
admet aussi une r\'esolution simpliciale, dont la filtration
(analogue \`a (0.1.2)) est homotopiquement triviale, voir [V2], [V4].
Le terme sup\'erieur non-trivial $\sigma_{m-1}\backslash \sigma_{m-2}$
de la filtration donne exactement le complexe des graphes connexes sur 
l'ensenble $M$, dont l'homologie correspond \`a la cohomologie
en degr\'e maximal de l'espace des applications 
injectives $M\hookrightarrow {\Bbb R}^d$. D'un autre c\^ot\'e son 
dual l'homologie en degr\'e maximal est d\'ecrite comme le 
sous-espace (d'une alg\`ebre de Lie libre), que l'on vient de d\'efinir, voir,
par exemple, [G], [Co].

La description de la cohomologie des complexes des graphes connexes
ainsi construite permet de d\'efinir le complexe dual au premier terme
de la suite spectrale auxiliaire.  L'homologie de ce complexe 
donne le premier terme
de la suite spectrale, qui est duale \`a la suite principale et qui calcule
l'homologie des espaces des n\oe uds. Dans ce travail je vais omettre
la plupart des d\'etails techniques (qui seront, d'ailleur, scrupuleusement 
\'etudi\'es dans ma th\`ese, voir [T]) et je donne tout de suite,
voir la section 1, la 
description de ce complexe, que l'on va appeler le complexe des
$*$-diagrammes de crochets ou bien le complexe des $B_*$-diagrammes 
et que l'on va d\'esigner par $CB_*D^{odd}(\Bbbk )$, $CB_*D^{even}(\Bbbk )$
({\it {\bf C}omplexe of {\bf B}racket $*$-{\bf D}iagrams}) pour $d$
impair et $d$ pair respectivement ($\Bbbk$ est un anneau commutatif de 
coefficients, $CB_*D^{odd(even)}(\Bbbk )\equiv 
CB_*D^{odd(even)}({\Bbb Z})\otimes{\Bbbk}$). 

Dans le discriminant on peut consid\'erer les strates, engendr\'ees
par les applications ${\Bbb R}\to {\Bbb R}^d$
avec seulement des auto-intersections
(on exclut les applications ayant des singularit\'es).
Les diagrammes dans $CB_*D^{odd(even)}(\Bbbk )$
 correspondant \`a tels strates sont appel\'es (simplement) {\it diag\-rammes
de crochets} ou {\it $B$-diagrammes}.
 L'espace engendr\'e par ces diagrammes
poss\`ede une structure de sous-complexe de $CB_*D^{odd(even)}(\Bbbk )$.
Le complexe ainsi obtenu est d\'esign\'e par   $CBD^{odd(even)}(\Bbbk )$
({\it {\bf C}omplexe of {\bf B}racket {\bf D}iagrams},
$CBD^{odd(even)}(\Bbbk )\equiv 
CBD^{odd(even)}({\Bbb Z})\otimes{\Bbbk}$).

Dans cette section on 
d\'ecrit \'egalement le moyen de simplifier les calculs
de l'homologie du complexe $CB_*D^{odd(even)}(\Bbbk )$. 
Ce complexe (qui simplifie les calculs et qui est homologiquement 
\'equivalent \`a $CB_*D^{odd(even)}(\Bbbk )$)
est un complexe-quotient de $CBD^{odd(even)}(\Bbbk )$; il est 
d\'esign\'e par $CB_0D^{odd(even)}(\Bbbk )$ ({\it {\bf C}omplexe 
of {\bf B}racket $0$-{\bf D}iagrams},
$CB_0D^{odd(even)}(\Bbbk )\equiv 
CB_0D^{odd(even)}({\Bbb Z})\otimes{\Bbbk}$).

Dans la section 2 on d\'efinit une structure d'alg\`ebres de Hopf 
diff\'erentielles (supercocommutatives) sur les complexes
$CB_*D^{odd(even)}(\Bbbk )$, $CBD^{odd(even)}(\Bbbk )$,
$CB_0D^{odd(even)}(\Bbbk )$. Les alg\`ebres de Hopf 
diff\'erentielles ainsi obtenues
sont d\'esign\'ees par $DHAB_*D^{odd(even)}(\Bbbk)$,
$DHABD^{odd(even)}(\Bbbk)$, $DHAB_0D^{odd(even)}(\Bbbk)$, respectivement.

En fait il arrive souvent, que la g\'eom\'etrie du discriminant
contienne des informations sur la structure comultiplicative
(et multiplicative, si l'espace de compl\'ement est un $H$-espace) dans 
l'homologie du compl\'ement.
Des conjectures sur le rapport en question dans le cas des
espaces de n\oe uds non-compacts sont formul\'ees dans la section 2.3
(Conjectures 2.3.5-6).

Si l'anneau principal $\Bbbk$ est un corps, alors l'homologie d'une alg\`ebre 
de Hopf diff\'erentielle (sur $\Bbbk$) forme une alg\`ebre de Hopf;
si $\Bbbk$ ne l'est pas, alors l'homologie correspondante
est \'etudi\'ee seulement comme une alg\`ebre sur $\Bbbk$.

L'espace de n\oe uds non-compacts dans ${\Bbb R}^d$, $d\ge 3$, est un 
$H$-espace, donc son homologie sur un corps forme une big\`ebre
(alg\`ebre de Hopf, pour $d\ge 4$). Si l'anneau de coefficients $\Bbbk$ 
n'est pas un corps, alors nous consid\'erons l'homologie (sur $\Bbbk$)
de ces espaces comme des alg\`ebres sur $\Bbbk$. On a les th\'eor\`emes 
suivants:

{\bf Th\'eor\`eme 2.3.7} [T]. {\it La big\`ebre de l'homologie sur un corps 
de caract\'eristique nulle de l'espace de n\oe uds non-compacts dans 
${\Bbb R}^d$, $d\ge 3$, est supercommutative.} $\Box$

{\bf Th\'eor\`eme 2.4.1.} {\it L'alg\`ebre (de Hopf) de l'homologie de
$DHABD^{odd(even)}(\Bbbk)$ est supercommutative pour n'importe quel
anneau commutatif $\Bbbk$.} $\Box$

{\bf Th\'eor\`eme 2.4.2.} {\it L'alg\`ebre (de Hopf) de l'homologie
de $DHAB_*D^{even}(\Bbbk)$ est supercommutative pour n'importe quel
anneau commutatif $\Bbbk$.} $\Box$

La m\'ethode de la d\'emonstration du Th\'eor\`eme 2.4.2 ne se
g\'en\'eralise pas dans le cas de $DHAB_*D^{odd}(\Bbbk)$.

Il est montr\'e \'egalement, voir la section 2.7, que l'inclusion
naturelle
$$
DHABD^{odd(even)}({\Bbbk})\hookrightarrow
DHAB_*D^{odd(even)}({\Bbbk})
\eqno(0.2.1)$$
pour le cas, o\`u $\Bbbk$ est un corps de caract\'eristique nulle,
induit une application injective en homologie. De plus,
pour $d$ pair le noyau est un id\'eal, engendr\'e par un seul g\'en\'erateur
primitif, pour $d$ impair -- par deux g\'en\'erateurs primitifs.
On en d\'eduit, que l'alg\`ebre de Hopf de l'homologie de 
$DHAB_*D^{odd}(\Bbbk)$ est supercommutative, si $\Bbbk$ est un corps de 
caract\'eristique nulle.

Le rapport avec la big\`ebre des diagrammes de cordes est donn\'e dans la 
section 2.8.

Dans les sections 3.1.1-2 on donne la construction (due \`a 
M.Gerstenhaber et A.Voronov [GV])
du complexe de Hochschild pour les op\'erades lin\'eaires 
munies d'un morphisme,
venant de l'op\'erade des alg\`ebres associatives. L'homologie de tels 
complexes poss\`ede une structure d'alg\`ebres de Gerstenhaber, i.e.
alg\`ebres supercommutatives munies d'un crochet impair de Lie 
compatible 
avec la multiplication.

Le r\'esultat principal de ce travail est, que les complexes
$CBD^{odd}(\Bbbk)$, $CBD^{even}(\Bbbk)$, $CB_*D^{even}(\Bbbk)$
sont les complexes de Hochschild normalis\'es pour les 
op\'erades d'alg\`ebres de Poisson, de Gerstenhaber, de Batalin-Vilkovissky,
respectivement (voir les Th\'eor\`emes 3.3.3, 3.3.6). Ce qui explique le fait
de la supercommutativit\'e dans l'homologie de $DHABD^{odd(even)}$ et de
$DHAB_*D^{even}$.

\vspace{4mm}

\noindent{\bf 0.3. Remerciments}

\vspace{2mm}

Premi\`erement je voudrais exprimer une reconnaissance profonde \`a mon 
directeur de recherches V.A.Vassiliev pour m'avoir pos\'e le 
probl\`eme, pour sa disponibilit\'e scientifique, pour son soutien 
et ses encouragements.

Je voudrais remercier sinc\`erement mes tuteurs --- Marc Chaperon (\`a Paris 7)
et  A.V.Chernavsky (\`a l'Universit\'e d'Etat de Moscou).

Je veux remercier M.Finkelberg et S.Loktev pour leurs tr\`es bons conseils 
math\'ematiques.

Et bien s\^ur, je suis tr\`es reconnaissant \`a P.Cartier, M.Kontsevich, 
M.D\'eza, D.Panov.

\vspace{4mm}

\noindent{\bf 0.4. Notations}

\vspace{2mm}

On note

$\bar{X}$ le compactifi\'e par un point d'un espace topologique
$X$;

$\tilde{H}_*(X),$ $\tilde{H}^*(X)$ la (co)homologie r\'eduite par rapport 
\`a un point d'un espace topologique $X$;

$d$ la dimension de l'espace consid\'er\'e ${\Bbb R}^d$;

$CB_*D^{odd}$, $CB_*D^{even}$ les complexes des $B_*$-diagrammes 
(= $*$-diagrammes de crochets) pour $d$ impair et $d$ pair respectivement;

$CBD^{odd}$, $CBD^{even}$ les complexes des $B$-diagrammes (= diagrammes de 
crochets);

$CB_0D^{odd}$, $CB_0D^{even}$ les complexes des $B_0$-diagrammes 
(= $0$-diagrammes de 
crochets);

Si l'on consid\`ere les cas de $d$ pair et impair simultan\'ement, on \'ecrit:
$CB_*D^{odd(even)}$, $CBD^{odd(even)}$, {\em etc}.

$DHAB_*D^{odd(even)}, DHABD^{odd(even)}$, $DHAB_0D^{odd(even)}$ les
alg\`ebres de Hopf diff\'erentielles des $B_*/B/B_0$-diagrammes.

Si l'on consid\`ere les cas des $*$-diagrammes et 
des diagrammes simultan\'ement, on \'ecrit:
($*$)-diagrammes, $B_{(*)}$-diagrammes, $CB_{(*)}D^{odd(even)}$, {\em etc.}

\newpage

\noindent{\bf 1. Les complexes des diagrammes de crochets, des 
$*$-diagrammes de crochets, des $0$-diagrammes de crochets
$CBD^{odd(even)}({\Bbbk})$,
$CB_*D^{odd(even)}({\Bbbk})$, $CB_0D^{odd(even)}({\Bbbk})$}           

\vspace{2mm}

\noindent{\bf 1.1. $(A,b)$-configurations}

\vspace{2mm}

Soit $A$ une collection finie non-ordonn\'ee de nombres entiers sup\'erieurs
ou \'egaux \`a 2, $A=(a_1,a_2,\dots,a_{\#A}),$ et soit 
$b$ un nombre entier non-n\'egatif. D\'efinissons 
$|A|:=a_1+\dots+a_{\#A}$. On va appeler $(A,b)$-{\it configuration} un ensemble
de $|A|$ points distincts de ${\Bbb R}$ d\'ecompos\'e en $\#A$ groupes
des cardinalit\'es $a_1,a_2,\dots,a_{\#A}$; et $b$ point diff\'erents
(dont une partie peut co\"\i ncider avec les $|A|$ points consid\'er\'es
ci-dessus). Ces $b$ points seront appel\'es {\it ast\'erisques}.
Pour la bri\'evit\'e les $(A,0)$-configurations seront appel\'ees simplement 
$A$-{\it configurations}.

Les $(A,b)$-configurations sont li\'ees \`a une stratification naturelle
dans $\sigma$ et $\Sigma$. Disons qu'une application 
$\Phi :{\Bbb R}\rightarrow{\Bbb R} ^d$, $\Phi\in \Sigma$,
 respecte une $(A,b)$-configuration, si
elle recolle les points dans chacun des groupes des cardinalit\'es 
$a_1,a_2, \dots,a_{\#A}$, et si sa d\'eriv\'e $\Phi '$ \'egale z\'ero
dans les derniers $b$ points de cette configuration.
Pour toute $(A,b)$-configuration l'ensemble des applications qui la respectent
est un sous-ensemble affine de ${\cal K}$ de codimension 
$d(|A|-\#A+b)$. Le nombre $|A|-\#A+b$ sera appel\'e la {\it complexit\'e} de 
la 
configuration. Deux $(A,b)$-configurations sont dites {\it \'equivalentes},
si elles peuvent \^etre obtenues l'une de l'autre par un hom\'eomorphisme
conservant l'orientation ${\Bbb
R}\rightarrow {\Bbb R}$.
Fixons une $(A,b)$-configuration $J$ de complexit\'e $i$ et avec $j$ points
g\'eom\'etriquement distincts sur ${\Bbb R}$. La strate form\'ee
 des applications ${\Bbb R}\to {\Bbb R}^d$ dans $\Sigma$,
qui respectent au moins une $(A,b)$-configuration $J'$ 
\'equivalente \`a $J$, peut \^etre param\'etris\'ee par l'espace
$S(J)$ de fibr\'e affine, dont la base est l'espace $E^j$
des $(A,b)$-configurations $J'$ \'equivalentes \`a $J$, et le fibr\'e
est ${\Bbb R}^{\omega -di}$. Notons que $E^j$ est contractible (etant une 
cellule ouverte de dimension $j$), alors cette fibration peut \^etre 
trivialis\'ee

$$
S(J)\simeq E^j\times {\Bbb R}^{\omega -di}.
\eqno(1.1.1)$$

{\bf Remarque 1.1.2.} La strate correspondante peut ne pas \^etre isomorphe
\`a $S(J)$ \`a cause de ses auto-intersections possibles. $\Box$

Permettons aux $(A,b)$-configurations (resp. $A$-configurations),
$A=(a_1,\dots,a_{\#A})$, 
d'avoir $a_i=1$ pour de certains $i\in\{1,\dots,\#A\}$.
De telles $(A,b)$-configurations seront appel\'ees $(A,b)$-configurations 
{\it g\'en\'eralis\'ees} (resp. $A$-configurations {\it g\'en\'eralis\'ees}).

Les $(A,b)$-configurations g\'en\'eralis\'ees avec au moins un 
$a_i=1$ (qui ne correspond \`a aucun ast\'erisque) n'ont pas 
d'interpr\'etation g\'eom\'etrique dans la stratification du discriminant,
mais elles nous seront tr\`es utiles 
pour certaines consid\'erations alg\'ebriques.

{\bf D\'efinition 1.1.3.} On va appeler {\it composante minimale} d'une 
$(A,b)$-configuration soit l'un de ses ast\'erisques qui ne co\"\i ncide
avec aucun des $|A|$ premiers points, soit l'un des $\#A$ groupes de points
avec tous les ast\'erisques qui y sont contenus. $\Box$

\vspace{4mm}

\noindent{\bf 1.2. Les espaces des diagrammes de crochets,
des $*$-diagrammes de crochets}

\vspace{2mm}

\noindent{\bf 1.2.1. Le cas o\`u $d$ est impair}

On fixe une $(A,b)$-configuration. Consid\'erons une superalg\`ebre de Lie
libre avec un crochet pair sur les g\'en\'erateurs impairs de type
$x_{t_\alpha}$, $x_{t^*_\beta}$, o\`u $t_\alpha$, 
$\alpha\in${\LARGE $\alpha$},
(resp. $t^*_\beta$, $\beta\in${\LARGE $\beta$}), sont tous les points sur 
${\Bbb R}$ de 
notre $(A,b)$-configuration, dans lesquels il n'y a pas d'ast\'erisques
(resp. il y en a)). On va prendre l'alg\`ebre 
sym\'etrique (dans le supersens) de 
l'espace de cette superalg\`ebre de Lie. Dans l'espace obtenu nous 
consid\`ererons un sous-espace engendr\'e par les produits de crochets, o\`u
chaque composante minimale de notre $(A,b)$-configuration (voir la
D\'efinition 1.1.3) est pr\'esent\'ee par un seul crochet tel, que 
tous les g\'en\'erateurs correspondant aux points de cette composante
y sont \'egalement pr\'esent\'es exactement une fois. De tels produits de crochets
seront appel\'es {\it $*$-diagrammes de crochets} ou bien 
{\it $B_*$-diagrammes}.

{\bf Exemple 1.2.1.} Pour la $(A,b)$-configuration de la Figure 1.2.2
on peut prendre le diagramme
$$
[[x_{t_1},x_{t_2^*}]x_{t_4}]\cdot[x_{t_3}x_{t_5}]\cdot
x_{t_6^*}\cdot x_{t_7^*}.~\Box
$$

\vspace{3mm}

\sevafig{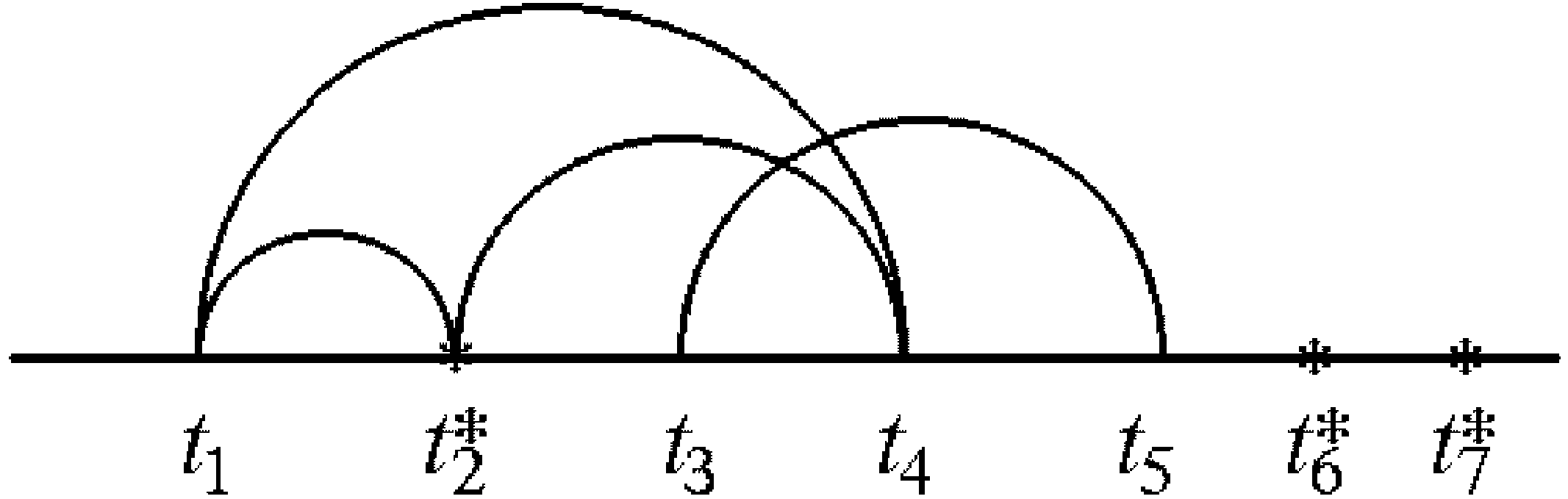}{10cm}

%\begin{figure}[htb]
%\center{\vbadness=10000\hbadness=10000\parbox[t]{50.1mm}{
%\hbox to 50.1mm{\vbox to 16.0mm{\special{em:graph 5.pcx}}}}}
%\parbox[b]{165mm}{
\centerline{(Figure 1.2.2)}
%}
%\end{figure}

\vspace{3mm}\vspace{3mm}

Les $B_*$-diagrammes correspondant aux $A$-configurations seront appel\'es 
simplement {\it diagrammes de crochets} ou {\it $B$-diagrammes}.

Deux ($*$)-diagrammes qui sont 
 obtenus l'un de l'autre par un hom\'eomorphisme
 ${\Bbb R}\to{\Bbb R}$, qui conserve l'orientation, sont consid\'er\'es comme
\'egaux.

{\bf Exemple 1.2.3.} Consid\'erons les $A$-configurations, dont toutes
les composantes minimales se composent de deux points (c'est exactement
les $(A,b)$-configurations qui donnent le nombre maximal ($=2i$) 
des points g\'eom\'etriquement distincts pour une complexit\'e fix\'ee $i$).
Les diagrammes de crochets, qui leur sont soumis, sont des produits de crochets
de type $[x_{t_1},x_{t_2}]$. De tels
diagrammes seront appel\'es {\it diagrammes de cordes}, parce que l'on peut 
les voir comme $2i$ points sur la droite, d\'ecompos\'es en $i$ couples
et joints par des cordes en dedans de toute couple. 
Notons que ni l'orientation de ces cordes ni leur ordonnancement
ne comptent, puisque
$[x_{t_1},x_{t_2}]=[x_{t_2},x_{t_1}]$ et que
ces crochets sont toujours pairs: $[x_{t_1},x_{t_2}]\cdot [x_{t_3},x_{t_4}]= 
[x_{t_3},x_{t_4}]\cdot[x_{t_1},x_{t_2}]$.   $\Box$

Les $B_*$-diagrammes ($B$-diagrammes), que l'on peut construire de la
m\^eme fa\c con pour les $(A,b)$-configurations (resp. $A$-configurations)
g\'en\'eralis\'ees, seront appel\'es $B_*$-diagrammes ($B$-diagrammes) 
{\it g\'en\'eralis\'es}.

{\bf Exemple 1.2.4.} Tout $B_*$-diagramme est un $B_*$-diagramme 
g\'en\'eralis\'e. $\Box$

{\bf Exemple 1.2.5.} $[x_{t_1},x_{t_2}]\cdot x_{t_3}$ est un $B$-diagramme
g\'en\'eralis\'e. $\Box$

\vspace{2mm}

\noindent{\bf 1.2.2. Le cas o\`u $d$ est pair}

Fixons une $(A,b)$-configuration et consid\'erons une superalg\`ebre de Lie
libre avec un crochet pair sur les g\'en\'erateurs pairs de type
$x_{t_\alpha}$ et les g\'en\'erateurs impairs de type $x_{t^*_\beta}$,
o\`u $t_{\alpha}$, $t^*_\beta$, $\alpha\in${\LARGE $\alpha$}, 
$\beta\in${\LARGE $\beta$}, sont
tous les points de notre $(A,b)$-configuration. Prenons l'alg\`ebre
 ext\'erieure (dans le supersens) de l'espace de cette alg\`ebre
de Lie. Dans l'espace obtenu on peut consi\-d\'e\-rer le sous-espace
lin\'eairement engendr\'e par les produits analogues (voir la section 1.2.1).
De tels produits de crochets seront \'egalement appel\'es 
{\it $B_*$-diagrammes}
({\it $*$-diagrammes de crochets}) et 
{\it $B$-diagrammes} ({\it diagrammes de crochets}).

{\bf Exemple 1.2.6.} Pour la $(A,b)$-configuration de la Figure 1.2.2 on peut
prendre le diagramme 
$$[[x_{t_1}x_{t_2^*}]x_{t_4}]\wedge[x_{t_3}x_{t_5}]\wedge
x_{t_6^*}\wedge x_{t_7^*}. ~\Box
$$

{\bf Exemple 1.2.7.} Pour $d$ pair les diagrammes, qui correspondent \`a 
l'Exemple 1.2.3, seront appel\'es {\it superdiagrammes de cordes}. 
Pour orienter tel diagramme il va importer l'orientation de ses cordes
(comme $[x_{t_1},x_{t_2}]=-
[x_{t_2},x_{t_1}])$ aussi bien que leur ordonnancement 
(comme $[x_{t_1},x_{t_2}]\wedge [x_{t_3},x_{t_4}]= 
-[x_{t_3},x_{t_4}]\wedge [x_{t_1},x_{t_2}]$)
. $\Box$

{\bf Remarque 1.2.8.} Il existe une autre mani\`ere de d\'efinir les 
$B_*$-diagrammes pour $d$ pair, qui donne la m\^eme chose. On prend 
la superalg\`ebre de Lie libre avec un crochet \underline{impair}
sur les g\'en\'erateurs \underline{impairs} $x_{t_\alpha}$, 
$\alpha\in${\LARGE $\alpha$}, 
et les g\'en\'erateurs \underline{pairs} $x_{t^*_\beta}$,
$\beta\in${\LARGE $\beta$}. 
Ensuite on consid\`ere l'alg\`ebre 
\underline{sym\'etrique} (dans le supersens) de l'espace de cette 
superalg\`ebre de Lie et on prend les produits (diagrammes) analogues. $\Box$

Par analogie (voir la section 1.2.1) on d\'efinit l'espace des
$B_*$-diagrammes ($B$-diagrammes) g\'en\'eralis\'es pour $d$ pair.

\vspace{4mm}

\noindent{\bf 1.3. La diff\'erentielle des complexes $CBD^{odd(even)}$, 
$CB_*D^{odd(even)}$ 
des ($*$)-diagrammes de crochets}

\vspace{2mm}

\noindent{\bf 1.3.0. Consid\'erations g\'en\'erales}

Les espaces des $B_*/B$-diagrammes ont une bigraduation naturelle.
La premi\`ere graduation est la complexit\'e $i$ des $(A,b)$-configurations
correspondantes. Le nombre $i$ indique l'homologie de quel terme 
$\sigma_i\backslash\sigma_{i-1}$ de la filtration (0.1.2) est calcul\'ee.
L'autre graduation est le nombre $j$ des points g\'eom\'etriquement 
diff\'erents des $(A,b)$-configurations. La diff\'erentielle sera de 
bidegr\'e (0,1). On va \'egalement associer un poids $p$ aux diagrammes:
$$
p:=i(d-1)-j.
$$

Selon les conjectures 0.1.3, 0.1.4 le complexe, que 
l'on va d\'efinir, devra calculer l'homologie (sur $\Bbb Q$)
de l'espace des n\oe uds. Le poids $p$ indique le degr\'e de l'homologie 
calcul\'ee.

Les formules (1.3.13), (1.3.7), (1.3.10),
(1.3.28), (1.3.31) donn\'ees ci-dessous
s'interpr\`ete g\'eom\`etriquement: elles d\'ecrivent
le cobord
des strates de type (1.1.1) parmi les strates de la m\^eme comlexit\'e $i$.

\vspace{2mm}

\noindent{\bf 1.3.1. Le cas o\`u $d$ est impair}

Soit $\gothg$ une superalg\`ebre de Lie munie d'un crochet pair:
$$
[.,.]:\gothg\otimes\gothg\rightarrow \gothg
\eqno(1.3.1)$$

Le crochet [.,.] admet une extension \`a l'alg\`ebre sym\'etrique
$S^*{\gothg} =\oplus_{i=0}^{+\infty}S^i{\gothg}$, que l'on appelle le {\it crochet
de Poisson}.

Soit $A,B\in S^*\gothg,$ o\`u $A=A_1\cdot A_2\dots A_k$, $B=B_1\cdot
B_2\dots B_\ell$, $A_i,B_j\in\gothg,$ $1\leq i\leq k,$ $1\leq j\leq
l,$ sont des \'el\'ements purs. On d\'efinit
$$
[A,B]:=\sum_{i,j}(-1)^{\lambda _{ij}}A_1\dots\hat{A}_i\dots
A_k\cdot[A_i,B_j]\cdot B_1\cdot\dots\cdot\hat{B}_j\cdot\dots\cdot B_\ell,
\eqno(1.3.2)$$

o\`u le chapeau sur un \'el\'ement signifie que ce dernier est omis;
$$
\lambda _{ij}=\tilde{A}_i\left(\sum_{p=i+1}^k
\tilde{A}_p\right)+\tilde{B}_j\left(\sum_{q=1}^{j-1}\tilde{B}_q\right).
\eqno(1.3.3)
$$

Le trait ondul\'e sur un \'el\'ement veut dire que l'on prend sa parit\'e.

Maintenant retournons \`a l'espace des ($*$)-diagrammes de crochets.
Soient $A$ et $B$ deux ($*$)-diagrammes de crochets, qui n'ont pas 
de points communs sur la droite ${\Bbb R}$.
Alors on peut aussi d\'efinir par la formule (1.3.2) un \'el\'ement
$[A,B]$ de l'espace des ($*$)-diagrammes de crochets.

{\bf D\'efinition 1.3.4.} On dit qu'une $(A,b)$-configuration $J$ peut \^etre
{\it ins\'er\'ee} en un point $t_0$ ou $t_0^*$ ({\it le signe ``$*$''
indique le fait, que ce point contient un ast\'erisque})
d'une $(A',b')$-configuration
$J'$, si $J$ n'a pas de points communs avec $J'$ \`a 
l'exception possible du point $t_0^{(*)}$. $\Box$

{\bf D\'efinition 1.3.5.} On dit qu'un diagramme peut \^etre {\it ins\'er\'e}
en un point $t_0^{(*)}$ d'un autre diagramme, si c'est vrai pour leurs 
$(A,b)$-configurations. $\Box$

Soient $A$ et $B$ deux ($*$)-diagrammes g\'en\'eralis\'es tels, que
$A$ puisse \^etre ins\'er\'e en un point $t_0$ (ou $t_0^*$) de $B$. On va 
d\'efinir un \'el\'ement 
$B|_{x_{t_0}=A}$ (resp. $B|_{x_{t^*_0}=A})$) de l'espace des 
$B_{(*)}$-diagrammes g\'en\'eralis\'es. Au signe pr\`es 
$B|_{x_{t^{(*)}_0}=A}$ est d\'efini par le remplacement formel de 
$x_{t^{(*)}_0}$ par $A$ (dans l'\'ecriture de $B$). Le signe est
d\'efini comme $(-1)^{(\tilde{A}-1)\times n}$, o\`u $n$ est le nombre des 
g\'en\'erateurs de type $x_{t_\alpha },x_{t^*_\beta }$ avant $x_{t^{(*)}_0}$
dans l'\'ecriture de $B$. En d'autres termes: On fait passer 
la composante minimale contenant $x_{t_0^{(*)}}$ \`a la premi\`ere place,
ensuite \`a l'aide des relations de superantisym\'etrie on met
$x_{t^{(*)}_0}$  \`a la premi\`ere place dans son crochet (et par 
cons\'equent dans l'\'ecriture de $B$); rempla\c cons $x_{t^{(*)}_0}$ par
$B$; et on fait toutes ces manipulations dans l'ordre inverse. Il est 
facile de 
voir que l'on obtient exactement le  signe, que l'on vient de d\'efinir.

{\bf Exemple 1.3.6.}
$$[x_{t_2}x_{t_3^*}]\cdot [x_{t_1^*}x_{t_0}]\bigl| _{x_{t_0}=[x_{t_4}x_{t_5}]
\cdot x_{t_6^*}}=(-1)^{(3-1)\cdot 3}[x_{t_2}x_{t_3^*}]\cdot
[x_{t_1^*},[x_{t_4}x_{t_5}]\cdot x_{t_6^*}]. ~\Box$$

Soient $A$ un $B_{(*)}$-diagramme non-g\'en\'eralis\'e, $t_\alpha$ 
l'un de ses points simples (=sans ast\'erisque), alors on d\'efinit 
$$ \partial_{t_\alpha
}A:=P\bigl(A|_{x_{t_\alpha }=x_{t_{\alpha -}}\cdot x_{t_{\alpha +}}}\bigr),
\eqno(1.3.7)$$
o\`u $P$ est la projection naturelle de l'espace des $B_{(*)}$-diagrammes
g\'en\'eralis\'es sur l'espace des  $B_{(*)}$-diagrammes, qui envoie en z\'ero 
tous les diagrammes ayant points simples isol\'es; on d\'efinit
$t_{\alpha\pm} par t_\alpha\pm\epsilon$ pour un tr\`es petit $\epsilon>0$.

{\bf Remarque 1.3.8.} La formule (1.3.7) peu \^etre pr\'esis\'ee:
$$\partial_{t_\alpha }A+(x_{t_{\alpha -}}-x_{t_{\alpha +}})\cdot
A=A|_{x_{t_\alpha }=x_{t_{\alpha -}}\cdot x_{t_{\alpha +}}}.~\Box
\eqno(1.3.9)$$

Soit $t^*_\beta $ l'un des points de $A$ (ayant un ast\'erisque), alors on 
d\'efinit
$$
\partial_{t_\beta ^*}A:=P\bigl(
A|_{x_{t^*_\beta }=x_{t_{\beta -}}\cdot x_{t^*_{\beta
+}}+x_{t^*_{\beta -}}\cdot x_{t_{\beta +}}+[x_{t_{\beta -}},x_{t_{\beta
+}}]}\bigr),
\eqno(1.3.10)$$
o\`u $P$ est la m\^eme projection; $t_{\beta\pm}^{(*)}:=t_\beta^*\pm\epsilon$.

{\bf Remarque 1.3.11.} La formule (1.3.10) peut \^etre pr\'ecis\'ee:
$$
\partial_{t^*_\beta }A+(x_{t_{\beta +}}-x_{t_{\beta +}})\cdot 
A=A|_{x_{t^*_\beta
}=x_{t_{\beta -}}\cdot x_{t^*_{\beta +}}+x_{t^*_{\beta -}}\cdot x_{t_{\beta
+}}+ [x_{t_{\beta -}},x_{t_{\beta +}}]}.~\Box
\eqno(1.3.12)$$

La diff\'erentielle $\partial$ sur l'espace des $B_*$-diagrammes sera d\'efinie
comme une somme des op\'erateurs $\partial_{t_\alpha}$ et
$ \partial_{t_\beta^*}$ par tous les points $t_\alpha$, $t_\beta^*$ de la
$(A,b)$-configuration correspondante:
$$
\partial =\sum_{\alpha\in
{\displaystyle \alpha }}\partial_{t_\alpha}+\sum_{\beta \in{\displaystyle
\beta }}\partial_{t_\beta^*}.
\eqno(1.3.13)$$

Il est facile de voir, que $\partial ^2=0$.

Le complexe des $B_*$-diagrammes (ainsi obtenus) est d\'esign\'e par 
$CB_*D^{odd}(\Bbbk )$, o\`u $\Bbbk$ est un anneau commutatif de coefficients.
L'espace des $B$-diagrammes d\'efinit un 
sous-complexe de $CB_*D^{odd}(\Bbbk )$,
qui sera d\'esign\'e par $CBD^{odd}(\Bbbk )$.

{\bf Remarque 1.3.14.}
$$\partial A= \left( \sum_{\alpha\in
{\displaystyle \alpha }}
A|_{x_{t_\alpha }=x_{t_{\alpha -}}\cdot x_{t_{\alpha +}}}\right)
+ \left( \sum_{\beta \in{\displaystyle
\beta }}A|_{x_{t^*_\beta
}=x_{t_{\beta -}}\cdot x_{t^*_{\beta +}}+x_{t^*_{\beta -}}\cdot x_{t_{\beta
+}}+ [x_{t_{\beta -}},x_{t_{\beta +}}]}\right) -
(x_{t_-}-x_{t_+})\cdot A,
\eqno(1.3.15)
$$
o\`u $t_+$ (resp. $t_-$) est un point sup\'erieur (resp. inf\'erieur) \`a 
tous les points de $A$. $\Box$

D\'efinissons \'egalement le complexe des $B_*$-diagrammes (resp. 
$B$-diagrammes) g\'en\'eralis\'es en prenant la diff\'erentielle $\partial$
selon la formule (1.3.15). On a une inclusion de $CB_*D^{odd}(\Bbbk )$
(resp. $CBD^{odd}(\Bbbk )$) dans le complexe ainsi d\'efini.

{\bf Affirmation 1.3.16.}[T] {\it Pour $d$ 
impair l'inclusion du complexe des $B_{(*)}$-diagrammes dans le complexe
des $B_{(*)}$-diagrammes g\'en\'eralis\'es induit un isomorphisme
en homologie.} $\Box$

{\bf D\'emonstration de 1.3.16:} Consid\'erons une filtration
d\'ecroissante dans le complexe des $B_{(*)}$-diagrammes g\'en\'eralis\'es
par le nombre des points isol\'es sans ast\'erisques.
Il est facile de voir, que ce complexe est une somme directe de 
$CB_{(*)}D^{odd}(\Bbbk )$ et du premier terme de cette filtration. On peut 
d\'emontrer, que le deuxi\`eme morceau
de cette somme est acyclique (voir [T]). $\Box$

{\bf Exemple 1.3.17.} Consid\'erons le complexe $CB_*D^{odd}(\Bbbk )$. C'est
clair que tous les dia\-grammes de cordes (voir l'Exemple 1.2.3) 
sont dans le noyau de la diff\'erentielle. Pour trouver les groupes
d'homologie dans les bigraduations $(i,2i)$ il faut quotienter
l'espace de tous les diagrammes de cordes par les relations, qui sont l'image
de la diff\'erentielle des diagrammes, dont toutes les composantes minimales,
sauf une seule, sont des cordes; celle qui ne l'est pas est un ast\'erisque
isol\'e ou bien un crochet sur 3 points sans ast\'erisques. 
Mentionnons que $\partial_{t_\beta^*}x_{t_\beta^*}=[x_{t_{\beta-
}},x_{t_{\beta+}}].$ Le cas d'un ast\'erisque isol\'e donne les dites 
{\it relations d'un terme} (voir la Figure 1.3.18).

\vspace{3mm}

\sevafig{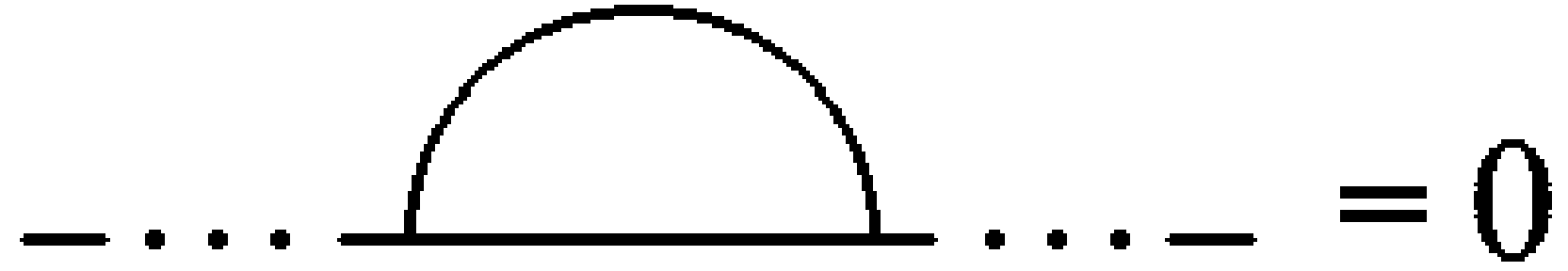}{6cm}
%\begin{figure}[htb]
%\center{\vbadness=10000\hbadness=10000\parbox[t]{33.8mm}{
%\hbox to 33.8mm{\vbox to 5.8mm{\special{em:graph 6.pcx}}}}}
%\parbox[b]{165mm}
\centerline{(Figure 1.3.18)}

%\end{figure}

\vspace{3mm}

Sur les lignes en pointill\'es il peut y avoir les bases d'autres cordes.
Autrement dit nous quotientons par les diagrammes de cordes ayant une 
corde joignant deux points voisins. Les calculs faciles de
$(\partial_{t_1}+\partial_{t_2}+\partial_{t_3})([[x_{t_1}x_{t_2}]x_{t_3
}])$ donnent les dites {\it relations de quatre termes} 
(voir la Figure 1.3.19).
$\Box$

\vspace{3mm}

\sevafig{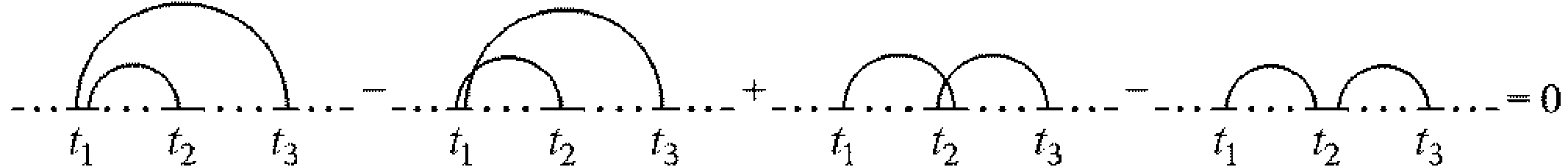}{12cm}

%\begin{figure}[htb]
%\center{\vbadness=10000\hbadness=10000\parbox[t]{109.7mm}{
%\hbox to 109.7mm{\vbox to 11.7mm{\special{em:graph 7.pcx}}}}}
%\parbox[b]{165mm}{
\centerline{(Figure 1.3.19)}

%\end{figure}

\vspace{3mm}

L'ordre des points $t_1$, $t_2$, $t_3$ sur la droite peut est arbitraire.
On met en pointill\'es les segments de la droite, 
o\`u il peut y avoir les bases
des autres cordes. (Les autres cordes sont toujours les m\^emes pour 
tous les quatre diagrammes.)

{\bf Exemple 1.3.20.} L'homologie de $CBD^{odd}(\Bbbk )$ dans les bidegr\'es
$(i,2i)$ est l' espace des diagrammes de cordes quotient\'e seulement
par les relations de quatre termes. $\Box$

\vspace{2mm}
\noindent{\bf 1.3.2. Le cas o\`u $d$ est pair}

Soit $\gothg$ une superalg\`ebre de Lie munie d'un crochet pair:
$$
[.,.]:\gothg\otimes\gothg\rightarrow \gothg.
\eqno(1.3.21)$$

Le crochet $[.,.]$ admet une extension sur l'alg\`ebre ext\'erieure
$\Lambda^*{\gothg}= \oplus_{i=0}^{+\infty}\Lambda^i{\gothg}$, que l'on appelle
le {\it crochet de Schouten}.

Soient $A=A_1\wedge\dots\wedge A_k\in\Lambda^k{\gothg}$,
$B=B_1\wedge \dots\wedge B_\ell\in\Lambda^\ell{\gothg}$. On d\'efinit
$$
[A,B]:=\sum_{i,j}(-1)^{\lambda_{ij}}A_1\wedge\dots\wedge \hat{A}_i
\wedge\dots\wedge A_k\wedge[A_i,B_j]\wedge
B_1\wedge\dots\wedge\hat{B}_j\wedge\dots\wedge B_\ell,
\eqno(1.3.22)$$

o\`u $\lambda_{ij}=(\tilde{A}_i+1)\left(\sum\limits_{p=i+1}^k(\tilde{A}_p
+1)\right)+(\tilde{B}_j+1)\left(\sum\limits_{q=1}^{j-1}(\tilde{B}_q+1)\right).
$

Il est facile de voir, que (1.3.22) d\'efinit correctement une application
$$
[\cdot,\cdot]:\Lambda^*{\gothg}\otimes\Lambda^*{\gothg}
\rightarrow\Lambda^*{\gothg}.
\eqno(1.3.23)$$

Pour tout mon\^ome de $\Lambda^*{\gothg}$ nous d\'efinissons sa parit\'e
comme la somme des parit\'es de ses facteurs plus le nombre des signes
de produit ext\'erieur. Par exemple, pour 
$A=A_1\wedge\dots\wedge A_k$
$$\tilde A=\tilde{A}_1+\dots+\tilde{A}_k+k-1.
\eqno(1.3.24)$$

Evidemment
$$
A\wedge B=(-1)^{(\tilde{A}+1)(\tilde{B}+1)}B\wedge A,
\eqno(1.3.25)$$

$$
[A,B]=-(-1)^{\tilde{A}\tilde{B}}[B,A].
\eqno(1.3.26)$$

Si $C=[A,B]$, alors $\tilde{C}=\tilde{A}+\tilde{B}$; si $C=A\wedge B$,
alors $\tilde{C}=\tilde{A}+\tilde{B}+1.$

Maintenant retournons \`a l'espace des ($*$)-diagrammes de crochets (pour le
cas o\`u $d$ est pair).

Soient $A=A_1\wedge\dots\wedge A_k,$ $B=B_1\wedge\dots\wedge
B_\ell$ deux ($*$)-diagrammes de crochets, qui n'ont pas de points communs
sur ${\Bbb R}$, $A_i,$ $B_j$ sont des crochets des g\'en\'erateurs pairs de 
type $x_{t_\alpha}$ et des g\'en\'erateurs impairs de type $x_{t_\beta^*}$.
D\'efinissons $[A,B]$ selon la formule (1.3.22). 

Notons que la parit\'e des diagrammes d\'efinie par (1.3.24) est exactement 
oppos\'ee \`a la parit\'e du poids des diagrammes (contrairement au cas o\`u 
$d$ est 
impair).

Soient $A$ et $B$ deux $B_{(*)}$-diagrammes g\'en\'eralis\'es, tels que 
$A$ peut \^etre ins\'er\'e dans un point $t_0^{(*)}$ de $B$. D\'efinissons
$B|_{x_{t_0^{(*)}}=A}$. Pour cela 
on remplace   $x_{t^{(*)}_0}$ par
$A$ dans l'\'ecriture de $B$, et on multiplie 
l' \'el\'ement ainsi obtenu  par 
$(-1)^{(\tilde{A}-\epsilon_0)\times (n_1+n_2)}$, 
o\`u $\epsilon_0$ \'egale $0$, s'il n'y a pas d'ast\'erisque
dans $t_0$, et \'egale $1$, s'il y a un ast\'erisque dans $t_0^*$,
$n_1$ (resp. $n_2$) est 
le nombre des signes de produit ext\'erieur  (resp. des
g\'en\'erateurs avec des ast\'erisques)
avant $x_{t_0^{(*)}}$ dans l'\'ecriture de $B$.

{\bf Exemple 1.3.27.}
$$[x_{t_2}x_{t_3^*}]\wedge [x_{t_1^*}x_{t_0}]\bigl| _{x_{t_0}=[x_{t_4}x_{t_5}]
\wedge x_{t_6^*}}=(-1)^{(2-0)\cdot (1+2)}[x_{t_2}x_{t_3^*}]\wedge
[x_{t_1^*},[x_{t_4}x_{t_5}]\wedge x_{t_6^*}]. ~\Box$$

Soient $A$ un $B_{(*)}$-diagramme, $t_\alpha$ l'un de ses points (simples),
alors d\'efinissons
 $$
\partial_{t_\alpha}A:=P\left(A|_{x_{t_\alpha}=x_{t_{\alpha-}}\wedge
x_{t_{\alpha+}}}\right),
\eqno(1.3.28)$$
o\`u $P$ est une projection analogue (\`a celle de (1.3.6)).

{\bf Remarque 1.3.29.} La formule (1.3.28) peut \^etre pr\'ecis\'ee:
$$
\partial_{t_\alpha}A+(x_{t_{\alpha-}}-x_{t_{\alpha+}})\wedge
A=A|_{x_{t_\alpha}=x_{t_{\alpha-}}\wedge x_{t_{\alpha+}}}.~ \Box
\eqno(1.3.30)$$

Soit $t_\beta^*$ l'un des points de $A$ (ayant un ast\'erisque), alors on 
d\'efinit
$$
\partial_{t^*_\beta} A:=P\left(A|_{x_{t^*_\beta}=x_{t_{\beta-}}\wedge
x_{t_{\beta+}^*}-x_{t_{\beta-}^*} \wedge x_{t_{\beta+}}-[x_{t_{\beta-
}},x_{t_{\beta+}}]}\right).~\Box
\eqno(1.3.31) $$

{\bf Remarque 1.3.32.} La formule (1.3.31) peut \^etre pr\'ecis\'ee:
$$
\partial_{t^*_\beta}A+(x_{t_{\beta-}}-x_{t_{\beta+}})\wedge
A=A|_{x_{t^*_\beta}=x_{t_{\beta-}}\wedge x_{t^*_{\beta+}}-
x_{t^*_{\beta-}}\wedge x_{t_{\beta+}}-[x_{t_{\beta-
}},x_{t_{\beta+}}]}.~\Box
\eqno(1.3.33)$$

D'une mani\`ere analogue au cas o\`u $d$ est impair on d\'efinit
l'op\'erateur $\partial$.
L'espace des $B_{(*)}$-diagrammes avec la diff\'erentielle $\partial$ est
appel\'e le {\it complexe des $B_{(*)}$-diagrammes} et est d\'esign\'e
par $CB_*D^{even}(\Bbbk )$, $CBD^{even}(\Bbbk )$. L'espace des 
$B_{(*)}$-diagrammes g\'en\'eralis\'es avec la diff\'erentielle $\partial$ est
appel\'e le {\it complexe des $B_{(*)}$-diagrammes g\'en\'eralis\'es.}

{\bf Affirmation 1.3.34.}[T] {\it Pour $d$ pair l'inclusion
naturelle de $CB_{(*)}D^{even}(\Bbbk )$ dans le
complexe des $B_{(*)}$-diagrammes g\'en\'eralis\'es
induit un isomorphisme en homologie.} $\Box$

{\bf D\'emonstration de 1.3.34:} analogue \`a celle de 
1.3.16. $\Box$

{\bf Remarque 1.3.35.} $$\partial A= \left( \sum_{\alpha\in
{\displaystyle \alpha }}
A|_{x_{t_\alpha }=x_{t_{\alpha -}}\wedge x_{t_{\alpha +}}}\right)
+ \left( \sum_{\beta \in{\displaystyle
\beta }}A|_{x_{t^*_\beta
}=x_{t_{\beta -}}\wedge x_{t^*_{\beta +}}-x_{t^*_{\beta -}}\wedge x_{t_{\beta
+}}- [x_{t_{\beta -}},x_{t_{\beta +}}]}\right) -
(x_{t_-}-x_{t_+})\wedge A.
$$
$\Box$

{\bf Exemple 1.3.36.} L'homologie des complexes $CB_*D^{even}(\Bbbk )$,
$CBD^{even}(\Bbbk )$ dans les bigra\-duations $(i,2i)$ est un espace
des superdiagrammes de cordes (voir l'Exemple 1.2.7) quotient\'e par les 
relations de quatre termes (voir la Figure 1.3.37) et 
(dans le cas de $CB_*D^{even}(\Bbbk )$ seulement) par les relations d'un terme.
$\Box$

\vspace{3mm}

\sevafig{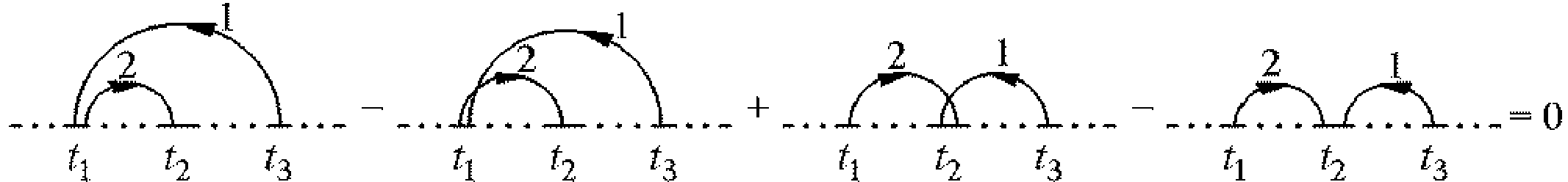}{12cm}

%\begin{figure}[htb]
%\center{\vbadness=10000\hbadness=10000\parbox[t]{112.4mm}{
%\hbox to 112.4mm{\vbox to 13.1mm{\special{em:graph 8.pcx}}}}}
%\parbox[b]{165mm}{
\centerline{(Figure 1.3.37)}

%\end{figure}

\vspace{3mm}

L'ordre des points $t_1$, $t_2$, $t_3$ peut \^etre arbitraire. Les nombres 
1 et 2 au-dessus des cordes d\'esignent l'ordre dans lequel on prend les 
cordes correspondantes. Les autres cordes et leur ordre sont les m\^emes
pour tous les quatre diagrammes. 

\vspace{4mm}

\noindent{\bf 1.4. Simplification des calculs de l'homologie de
$CB_*D^{odd(even)}(\Bbbk )$. 
Les complexes $CB_0D^{odd(even)}(\Bbbk )$
des $B_0$-diagrammes}

\vspace{2mm}

Consid\'erons l'espace des $B$-diagrammes (pour $d$ pair ou impair).
On va le quotienter par les relations de ``supercommutativit\'e
voisine'', plus pr\'ecisement aux relations habituelles de 
superantisym\'etrie et celles  de Jacobi on ajoute la 
supercommutabilit\'e des g\'en\'erateurs correspondant aux points
voisins sur la droite ${\Bbb R}$ dans 
les diagrammes. L'espace ainsi obtenu sera appel\'e
{\it l'espace des $B_0$-diagrammes} (ou bien
{\it l'espace des $0$-diagrammes de crochets}). L'espace des
 $B_0$-diagrammes
poss\`ede une structure de complexe-quotient de $CBD^{odd(even)}(\Bbbk )$
et de
$CB_*D^{odd(even)}(\Bbbk )$ (dans le deuxi\`eme cas on envoie en z\'ero
tous les diagrammes ayant des ast\'erisques), que l'on appelle le 
{\it complexe des $B_0$-diagrammes} et qui est d\'esign\'e par 
$CB_0D^{odd(even)}(\Bbbk )$.

On a un diagramme commutatif des morphismes de complexes:
$$
\begin{array}{c}
CB_0D^{odd(even)}(\Bbbk)  \stackrel{\displaystyle p_1}{\leftarrow }
CBD^{odd(even)}(\Bbbk)\\ p_2\nwarrow \quad \swarrow i_1\\
 CB_*D^{odd(even)}(\Bbbk)
\end{array}
\eqno(1.4.1)$$

o\`u $p_1$, $p_2$ sont projectifs; $i_1$ est injectif.

{\bf Th\'eor\`eme 1.4.2.}[T] {\it L'espace des $B_0$-diagrammes 
(sur n'importe quel anneau commutatif $\Bbbk$ de coefficients) est un
$\Bbbk$-module libre (autrement dit, ce, que l'on quotiente, ne donne pas 
de torsions). La surjection $p_2$ de (1.4.1) induit un isomorphisme 
en homologie.} $\Box$

{\bf Id\'ee de  la d\'emonstration:} Il faut consid\'erer la filtration par 
le nombre des composantes minimales dans le complexe dual \`a
$CB_*D^{odd(even)}(\Bbbk )$ et la suite spectrale associ\'ee
\`a cette filtration. Cette suite d\'eg\'en\`ere au deuxi\`eme terme, 
parce que son premier terme se trouve dans une seule ligne, qui r\'epond 
aux diagrammes sans ast\'erisques. Pour les d\'etails voir [T]. $\Box$

\vspace{4mm}

\noindent{\bf 2. Sur six alg\`ebres de Hopf diff\'erentielles li\'ees
aux discriminants des espaces de n\oe uds non-compacts}

\vspace{2mm}

\noindent{\bf 2.1. Sur les alg\`ebres de Hopf (diff\'erentielles)}

\vspace{2mm}

Soit $(B,\mu,\iota,\Delta,\epsilon)$ une big\`ebre, o\`u $B$ est un
module ${\Bbb Z} _2$-gradu\'e sur un anneau commutatif $\Bbbk$;
$\mu,\iota,\Delta,\epsilon$ sont respectivement 
les application de la multiplication, de l'unit\'e, de la comultiplication,
de la counit\'e.

Sur l'espace ${\Bbb Z} _2$-gradu\'e $Mor(B,B)$ de toutes les applications 
lin\'eaires de $B$ dans $B$ on a une structure d'alg\`ebre associative 
en d\'efinissant la multiplication $\star$ par la formule
 $$
f\star g:=\mu\circ(f\otimes g)\circ\Delta.
\eqno(2.1.1)$$
L'unit\'e de cette alg\`ebre est
$$
\1 =\iota\circ\epsilon.
\eqno(2.1.2)$$

Soit $S\in Mor(B,B)$ tel, que
$$
S\star id=id\star S=\1,
\eqno(2.1.3)$$
alors $S$ est appel\'e {\it antipode}. Evidemment, s'il existe, il est unique.
Tout tel sextuplet  $(B,\mu,\iota,\Delta,\epsilon,S)$ est appel\'e
{\it alg\`ebre de Hopf}.

Rappelons  qu'une  {\it alg\`ebre de Hopf diff\'erentielle}
$(H,\mu,\iota,\Delta,\epsilon,S,\partial )$ est une alg\`ebre de Hopf
 $(H,\mu,\iota,\Delta,\epsilon,S)$ sur laquelle agit une application
impaire $\partial\in Mor(B,B)$, v\'erifiant $\partial^2=0$ et
compatible avec les applications de structure $\mu$, $\iota$, 
$\Delta$, $\epsilon$, $S$.

De la m\^eme fa\c con on d\'efinit big\`ebre, alg\`ebre, coalg\`ebre 
diff\'erentielles.

L'homologie (par rapport \`a la diff\'erentielle $\partial$) d'une alg\`ebre
diff\'erentielle sur n'importe quel anneau commutatif $\Bbbk$ forme 
\'egalement une alg\`ebre. L'homologie d'une coalg\`ebre diff\'erentielle
peut ne pas \^etre une coalg\`ebre, si l'anneau principal
n'est pas un corps. De la m\^eme fa\c con l'homologie d'une big\`ebre
(alg\`ebre de Hopf) diff\'erentielle sera consid\'er\'ee comme une big\`ebre
(alg\`ebre de Hopf), si l'anneau principal est un corps, et
simplement comme une alg\`ebre, s'il ne l'est pas.

{\bf D\'efinition 2.1.4.} Une big\`ebre non-n\'egativement $\Bbb Z$-gradu\'ee 
est
dite {\it connexe}, si son espace de degr\'e 0 est engendr\'e par un 
seul \'el\'ement non-nul (ce qui revient \`a dire qu'il est isomorphe \`a
un anneau-quotient de l'anneau principal $\Bbbk$). $\Box$

Pour toute big\`ebre connexe $B$ les op\'erateurs de la forme
$$
\sum_{k=0}^{+\infty}a_k(id-\1)^{\star k},
\eqno(2.1.5)$$
(o\`u $(id-\1)^{\star k}=\underbrace{(id-\1)\star(id-\1)
\star\dots\star(id-\1 )}_{k};$
$a_k\in\Bbbk$, $k=0,1,2,\dots$) sont correctement d\'efinis, parce que
pour tout \'el\'ement $x\in B$ 
on a $(id-\1 )^{\star k}x=0$, si $k$ est sup\'erieur
au degr\'e de $x$.

Soit $f(t)=\sum\limits_{k=0}^{+\infty}a_k(t-1)^k$ la fonction g\'en\'eratrice
d'une suite $\{ a_k\}$. D\'efinissons
$$
f_\star(id):=\sum_{k=0}^{+\infty}a_k(id-\1)^{\star k}.
\eqno(2.1.6)$$

{\bf Lemme 2.1.7.} {\it Soit $(B,\mu ,\iota,\Delta ,\epsilon ,\partial)$ une 
big\`ebre diff\'erentielle, soient  $f,g\in Mor(B,B)$
des op\'erateurs compatibles avec la diff\'erentielle $\partial$ 
(supercommutant avec elle), alors $f\star g$ est \'egalement compatible avec
$\partial$.} $\Box$

{\bf D\'emonstration du Lemme 2.1.7:} C'est une cons\'equence
du fait que $\mu ,\Delta ,f,g$
sont compatibles avec $\partial$. $\Box$

{\bf Cons\'equence du Lemme 2.1.7.}  Tout op\'erateur de la forme (2.1.5)
dans une big\`ebre diff\'erentielle connexe est toujours compatible avec la 
diff\'erentielle. $\Box$

{\bf Th\'eor\`eme 2.1.8.} {\it Toute big\`ebre (diff\'erentielle) connexe sur 
n'importe quel anneau commutatif $\Bbbk$ est une alg\`ebre de Hopf 
(diff\'erentielle). L'antipode $S$ est donn\'e par la formule suivante:}
$$
S=\sum_{k=0}^{+\infty}(-1)^k(id-\1)^{\star k}.~\Box
\eqno(2.1.9)$$

{\bf D\'emonstration du Th\'eor\`eme 2.1.8:} Evident. $\Box$

{\bf Th\'eor\`eme 2.1.10.}[P] {\it Toute big\`ebre diff\'erentielle 
supercocommutative connexe sur un corps de caract\'eristique nulle 
est isomorphe
comme alg\`ebre de Hopf diff\'erentielle
 \`a l'alg\`ebre enveloppante
de
la superalg\`ebre de Lie diff\'erentielle de ses \'el\'ements primitifs.
L'op\'erateur
$$ P_1=\log _\star id
\eqno(2.1.11)
$$
est une projection sur l'espace des \'el\'ements primitifs. L'homologie
de telle alg\`ebre de Hopf diff\'erentielle est l'alg\`ebre de Hopf 
enveloppante de l'homologie 
de la superalg\`ebre de Lie diff\'erentielle en question.} $\Box$

\vspace{4mm}

\noindent{\bf 2.2. Six alg\`ebres de Hopf diff\'erentielles $DHAB_*D^{odd},$
$DHAB_*D^{even}$, $DHAB_0D^{odd}$, $DHAB_0D^{even}$, $DHABD^{odd}$,
$DHABD^{even}$ des ($*$/$0$)-diagrammes de crochets}

\vspace{2mm}

Consid\'erons les complexes $CB_*D^{odd(even)}(\Bbbk),$ 
$CB_0D^{odd(even)}(\Bbbk)$,
$CBD^{odd(even)}(\Bbbk)$
des ($*$/$0$)-diagrammes de crochets. Ces complexes sont bigradu\'es
par la complexit\'e $i$ et par le nombre des points g\'eom\'etriquement
diff\'erents. La diff\'erentielle y est de bidegr\'e (0,1).

On va construire les op\'erations $\mu ,\iota,\Delta ,\epsilon ,S$ de
multiplication, unit\'e, comultiplication, counit\'e, antipode, qui 
 vont d\'efinir avec
$\partial$ une structure d'alg\`ebre de Hopf 
diff\'erentielle bigradu\'ee
sur ces espaces. Ces alg\`ebres de Hopf diff\'erentielles seront appel\'ees
respectivement $DHAB_*D^{odd(even)}(\Bbbk)$, $DHAB_0D^{odd(even)}(\Bbbk)$,
$DHABD^{odd(even)}(\Bbbk)$ (de l'anglais ``{\bf D}ifferential 
{\bf H}opf {\bf A}lgebra of
{\bf B}racket {\bf D}iagrams'').

D\'efinissons $\iota:\Bbbk\rightarrow DHAB_xD^{odd(even)}(\Bbbk )$ en
prenant $\iota(1)$ \'egale le diagramme trivial (\underline{ici 
et plus loin
$x$ signifie soit $*$, soit 0, soit $\emptyset$}).

D\'efinissons $\epsilon$ comme z\'ero sur tous les diagrammes
sauf celui qui est trivial. $\epsilon$ du diagramme trivial est d\'efini 
comme 1.

Notos que $\iota,\epsilon $ sont les isomorphismes entre l'anneau
$\Bbbk$ et le module de bidegr\'e (0,0). Donc 
$\iota,\epsilon $ respectent la bigraduation. Evidemment ils respectent
\'egalement la diff\'erentielle $\partial$, puisque la diff\'erentielle
$\partial$ du diagramme
trivial est nulle.

D\'efinissons la multiplication $\mu$. Le produit de deux diagrammes
$D_1$ et $D_2$ est d\'efini comme un diagramme $D_1*D_2=\mu (D_1\otimes D_2)$
sur la droite ${\Bbb R}$, que l'on obtient en collant le diagramme
$D_2$ \`a droite du diagramme $D_1$ (voir la Figure 2.2.1).
%\vspace{3mm}

%\begin{figure}[htb]
\hspace*{25mm}
\PSbox{pic1.pstex}{10mm}{15mm} 
\begin{picture}(10,10)
\put(-12,12){$D_1$}
\put(50,12){$D_2$}
\put(130,12){$D_1$}
\put(194,12){$D_2$}
\end{picture}

{\centerline{(Figure 2.2.1)}
%\end{figure}

%\begin{figure}[htb]
%\center{\vbadness=10000\hbadness=10000\parbox[t]{61.0mm}{
%\hbox to 61.0mm{\vbox to 6.6mm{\special{em:graph 9.pcx}}}}}
%\parbox[b]{165mm}
%{\centerline{(Figure 2.2.1)}
%}
%\end{figure}

\vspace{3mm}

De la m\^eme fa\c{c}on on obtient l'\'ecriture en crochets de $D_1*D_2$
en ajoutant \`a droite de l'\'ecriture de $D_1$ celle de $D_2$.

{\bf Exemple 2.2.2.} $\bigl( [x_{t_1},x_{t_2}]\wedge [x_{t_3},x_{t_4^*}]
\bigr) \star [x_{t_1},[x_{t_2},x_{t_3}]]=[x_{t_1},x_{t_2}]\wedge 
[x_{t_3},x_{t_4^*}]
\wedge [x_{t_1+N},[x_{t_2+N},x_{t_3+N}]]$
pour un tr\`es grand $N> 0$. $\Box$

D\'efinissons la comultiplication $\Delta$. Soit $D$ un ($*$/$0$)-diagramme
de crochets, soit $V=V(D)$ l'ensemble de ses composantes minimales.
Tout sous-ensemble $V_1\subset V$ de composantes minimales du diagramme
$D$ d\'efinit un diagramme $(V_1)$, qui ne se compose que des composantes 
minimales de $V_1$ (on efface tous les crochets dans l'\'ecriture de $D=(V)$
qui correspondent aux composantes minimales, qui ne sont pas de $V_1$).

On d\'efinit $\Delta $ en prenant
$$
\Delta((V))=\sum_{V=V_1\sqcup V_2}\pm(V_1)\otimes(V_2).
\eqno(2.2.3)$$
Le signe $\pm$ dans (2.2.3) peut \^etre pr\'ecis\'e facilement:
il
appara\^\i t, quand on fait passer \`a droite
les composantes minimales de $V_2$ dans l'\'ecriture de $D$.

Il est facile de voir, que la comultiplication est supercocommutative.

On peut d\'emontrer, que les op\'erations $\mu ,\iota,\Delta ,\epsilon
,\partial$ font de $CB_xD^{odd(even)}$ une big\`ebre diff\'erentielle connexe
(voir la D\'efinition 2.1.4) $\Bbb Z$-bigradu\'ee.
Cela veut dire (voir le Th\'eor\`eme 2.1.8),
qu'elle est en fait une alg\`ebre de Hopf 
diff\'erentielle (il y existe un antipode $S$).

{\bf D\'efinition 2.2.4.} Un intervalle de ${\Bbb R}$ 
joignant deux points voisins d'une $(A,b)$-configuration est 
dit {\it s\'eparant} pour elle,
si toute composante minimale de cette $(A,b)$-configuration 
est soit strictement \`a gauche, soit strictement \`a droite par 
rapport \`a cet intervalle-l\`a. $\Box$

{\bf D\'efinition 2.2.5.} En \'eliminant tous les intervalles s\'eparants 
dans la droite ${\Bbb R}$, on obtient une d\'ecomposition de notre 
$(A,b)$-configuration en {\it composantes connexes}, qui r\'epondent
aux secteurs connexes de la droite sans ces intevalles. $\Box$

Une composante connexe peut contenir plusieures composantes minimales.

{\bf Exemple 2.2.6.} La $(A,b)$-configuration de la Figure 1.2.2 a 4 
composantes minimales et 3 composantes connexes. $\Box$

Consid\'erons l'espace ${\cal P}$ engendr\'e par les diagrammes de 
$DHAB_xD^{odd(even)}(\Bbbk )$ n'ayant qu'une seule composante 
connexe. Il est facile de voir que  l'alg\`ebre
$DHAB_xD^{odd(even)}(\Bbbk )$ 
est l'alg\`ebre libre tensorielle
$T^*({\cal P})$ de l'espace $\cal P$. Soit $\Bbbk$ un corps de 
caract\'eristique nulle, alors on a une projection
$P_1=\log_\star id$ sur l'espace des \'el\'ements primitifs
(voir le Th\'eor\`eme 2.1.10).
Il est clair que ${\rm ker}\,P_1|_{\cal P}=0.$ On en d\'eduit que l'alg\`ebre
$DHAB_xD^{odd(even)}(\Bbbk)$ est
isomorphe \`a $T^*(P_1({\cal P}))$.
Cela veut dire que l'alg\`ebre de Hopf
$DHAB_xD^{odd(even)}(\Bbbk)$  
est isomorphe \`a l'alg\`ebre de Hopf libre
d'un certain ensemble infini de g\'en\'erateurs pairs et impairs primitifs
$\Bbb Z$-bigradu\'es, ou bien encore
$DHAB_xD^{odd(even)}(\Bbbk)$ est l'alg\`ebre enveloppante d'une superalg\`ebre
de Lie libre d'un en\-semble infini de g\'en\'erateurs pairs et impairs
$\Bbb Z$-bigradu\'es. D'apr\`es le Th\'eor\`eme 2.1.10 l'alg\`ebre de Hopf 
de l'homologie de $DHAB_xD^{odd(even)}(\Bbbk)$ doit \^etre l'alg\`ebre 
enveloppante d'une superalg\`ebre de Lie bigradu\'ee. Il se trouve
que les quatre superalg\`ebres de Lie ainsi obtenues sont toutes 
supercommutatives (ont un crochet trivial).

Notons que les morphismes $i_1$, $p_1$, $p_2$ du diagramme (1.4.1)
sont aussi des morphismes des alg\`ebres de Hopf diff\'erentielles. On a
donc le diagramme commutatif suivant:
$$
\begin{array}{c}
DHAB_0D^{odd(even)}(\Bbbk)\stackrel{p_1}{\longleftarrow}DHABD^{odd(even)}(\Bbbk)\\
p_2\nwarrow\quad\swarrow i_1\\
DHAB_*D^{odd(even)}(\Bbbk)\\
\end{array}
\eqno(2.2.7)$$

L'application $p_2$ induit un isomorphisme en homolgie.

\vspace{2mm}

\noindent{\bf 2.3. Conjectures sur la multiplication et la comultiplication}
{\bf dans la (co)homologie des espaces de n\oe uds non-compacts}

\vspace{2mm}

Un espace topologique $X$ sera appel\'e {\it $H$-espace}, s'il est muni
d'une multiplication homotopiquement associative
$$
m:X\times X\rightarrow X,
\eqno(2.3.1)$$
et s'il a un \'el\'ement $e\in X$, qui est une unit\'e 
homotopique par rapport \`a $m$. L'existence d'une application 
homotopiquement inverse n'est pas d\'emand\'ee.

Il est facile de voir, que l'homologie et la cohomologie sur un corps
de tout $H$-espace forment des big\`ebres respectivement supercocommutative et
supercommutative, qui sont duales l'une \`a l'autre. Si notre $H$-espace
est connexe par arcs ou bien si l'on peut d\'efinir pour lui une application
homotopiquement inverse, alors ces big\`ebres sont en fait des alg\`ebres 
de Hopf (dans le premier cas c'est grace au Th\'eor\`eme 2.1.8).

L'homologie et la cohomologie sur un anneau commutatif $\Bbbk$
(qui n'est pas un corps) des $H$-espaces seront consid\'er\'ees comme des 
alg\`ebres sur $\Bbbk$.

Les espaces de n\oe uds non-compacts dans ${\Bbb R}^d$, $d\ge 3$, sont des 
$H$-espaces: on prend la composition des n\oe uds comme multiplication et
 le n\oe ud trivial comme unit\'e. Alors leur (co)homologie sur un corps est
toujours une big\`ebre (alg\`ebre de Hopf, pour $d\ge 4$).

La filtration (0.1.2) dans le discriminant r\'esolu $\sigma$ d\'efinit une 
filtration dans l'homologie des espaces 
$\bar{\sigma}\simeq\bar{\Sigma}:$
$$
\begin{array}{ccccccc}
\tilde{H}_*^{(0)}(\bar{\sigma},\Bbbk) & \subset &\tilde{H}_*^{(1)}
(\bar{\sigma}, \Bbbk) &
\subset & \tilde{H}_*^{(2)}(\bar{\sigma},\Bbbk) & \subset & \dots\\
\textstyle |\wr && |\wr && |\wr &&\\
\tilde{H}_*^{(0)}(\bar{\Sigma},\Bbbk) &\subset & 
\tilde{H}_*^{(1)}(\bar{\Sigma}, \Bbbk)
&\subset  & \tilde{H}_*^{(2)}(\bar{\Sigma},\Bbbk) & \subset &\dots .  
\end{array}
\eqno(2.3.2)$$

Grace \`a l'isomorphisme d'Alexander cela donne une filtration 
croissante dans la cohomologie des espaces de n\oe uds
$$H_{(0)}^*({\cal K}\backslash\Sigma,\Bbbk)\subset H^*_{(1)}({\cal
K}\backslash\Sigma, \Bbbk)\subset H_{(2)}^*({\cal
K}\backslash\Sigma,\Bbbk)\subset \dots .
\eqno(2.3.3)$$

De la m\^eme fa\c con on obtient une filtration d\'ecroissante
dans l'homologie
$$
H_*^{(0)} ({\cal K}\backslash\Sigma,\Bbbk)\supset H_*^{(1)} ({\cal
K}\backslash\Sigma,\Bbbk)\supset H_*^{(2)} ({\cal
K}\backslash\Sigma,\Bbbk)\supset\dots .
\eqno(2.3.4)$$
Pour $d\ge 4$ les filtraions (2.3.3-4) sont toujours finies pour tout degr\'e
$*$; si $d=3$ la filtration (2.3.3) n'\'epuise pas toute la cohomologie 
de l'espace des n\oe uds. La conjecture de la compl\'etude
de la ``(co)homologie de type fini'' est en fait la condition de la 
convergence de la filtration (2.3.4) vers z\'ero. Ce probl\`eme
n'est toujours pas r\'esolu m\^eme pour les invariants -- la cohomologie 
de degr\'e z\'ero.

{\bf Conjecture 2.3.5.} {\it La multiplication (et la comultiplication,
si l'anneau principal $\Bbbk$ est un corps) dans la cohomologie
et l'homologie respectent les filtrations (2.3.3) et (2.3.4)
respectivement.} $\Box$

Modulo cette conjecture on d\'efinit une structure d'alg\`ebre (big\`ebre)
sur les qoutients gradu\'es (par rapport aux filtration (2.3.3-4)).

De la formule (2.1.9) et de l'hypoth\`ese que la multiplication
et la comultiplication respectent la filtration on d\'eduit, que 
l'antipode $S$ la respecte \'egalement, et donc il est aussi induit
sur le quotient gradu\'e.

La suite spectral principale (voir la section 0.1) et la suite, qui lui est
 duale, 
convergent, evidemment (comme $\Bbbk$-modules) vers ces quotients gradu\'es.
Dans la section pr\'ecedente on a construit une structure d'alg\`ebres de Hopf
diff\'erentielles sur les complexes $CB_*D^{odd(even)}(\Bbbk )$, qui
sont isomorphes au premier terme de la suite spectrale auxiliaire
(voir la section 0.1). De la m\^eme fa\c con on peut
d\'efinir une structure d'alg\`ebres de Hopf diff\'erentielles
sur le z\'eroi\`eme terme de la suite principale et sur 
celui de la suite duale.
Nous n'allons pas d\'ecrire les z\'eroi\`emes termes en question, cependant
nous formulerons la conjecture suivante.

{\bf Conjecture 2.3.6.} {\it La structure d'alg\`ebre (alg\`ebre de Hopf, 
si l'anneau principal est un corps) diff\'erentielle sur le z\'eroi\`eme
terme de la suite spectrale principale (ou de la suite duale) est
 induite sur tous les termes $E_{*,*}^1,
E_{*,*}^2, \dots (E_1^{*,*},E_2^{*,*},\dots)$, o\`u pour une diff\'erentielle
sur le $i$-\`eme terme on prend la $i$-\`eme diff\'erentielle de 
la suite spectrale. La structure d'alg\`ebre (de Hopf) sur 
$E_{*,*}^\infty$ (resp.
$E_\infty^{*,*}$) ainsi obtenue co\"\i ncide avec celle (du quotient gradu\'e),
que l'on a conjectur\'ee dans (2.3.5).} $\Box$

Si $\Bbbk$ est un corps de caract\'eristique nulle,
alors selon la conjecture 0.1.3 
la suite spectrale principale et sa duale d\'eg\'en\`erent au premier terme.
On en d\'eduit (modulo les Conjecture 2.3.5-6),
 que l'alg\`ebre de Hopf de l'homologie de
$DHAB_*D^{odd(even)}(\Bbbk )$ est exactement le quotient gradu\'e, qui nous 
int\'eresse. Si l'on suppose aussi (la Conjecture 0.1.4 de Vassiliev),
que la filtration (0.1.2) est homotopiquement trivial, alors ce sera vrai pour
tout anneau commutatif $\Bbbk$.

Dans [T] on d\'emontrera le th\'eor\`eme suivant.

{\bf Th\'eor\`eme 2.3.7.} [T] {\it La big\`ebre de l'homologie sur $\Bbb Q$
des espaces de n\oe uds non-compacts dans ${\Bbb R}^d$, $d\ge 3$,
est toujours supercommutative.} $\Box$

{\bf Id\'ee de la d\'emonstration:} Si l'on a deux cycles, qui representent
certaines classes de l'homologie $H_*({\cal K}\backslash\Sigma ,\Bbbk )$,
et l'on consid\`ere leur produit, voir la Figure 2.3.8, alors on

\hspace*{10mm}
\PSbox{knot0.pstex}{15mm}{48mm} 
\begin{picture}(10,10)
\end{picture}

\vspace{3mm}\vspace{1mm}

\parbox[b]{165mm}{
\centerline{(Figure 2.3.8)}
}
\vspace{3mm}

\noindent va comprimer celui, qui est 
\`a droite, jusqu'\`a ce qu'il soit tr\`es petit, et ensuite on le fait passer 
(\`a l'aide de la connexion adiabatique) par le cycle gauche. Le cycle
ainsi obtenu ne sera point le produit du deuxi\`eme cycle sur le premier,
parce que la connexion adiabatique donne une perturbation ---
une rotation de $SO(d-1)$, qui d\'epend des n\oe uds du premier cycle.
Il se trouve (voir [T]), que cette perturbation ne vas pas importer,
si $\Bbbk$ est un corps de caract\'eristique nulle. $\Box$

Il suit de ce th\'eor\`eme (toujours modulo les Conjectures 0.1.3, 2.3.5,
2.3.6),
que l'alg\`ebre de Hopf de l'homologie de  $DHAB_*D^{odd(even)}({\Bbb Q} )$
est supercommutative. Cette affirmation sera d\'emontr\'ee ind\'ependamment
(voir les sections 2.4, 2.7).

{\bf Remarque 2.3.9.} Du Th\'eor\`eme 2.3.7 (modulo la Conjecture 0.1.3)
on d\'eduit, que la big\`ebre de l'homologie
sur $\Bbb Q$ de l'espace de n\oe uds non-compacts dans ${\Bbb R}^d$,
$d\ge 4$, est isomorphe \`a la big\`ebre de l'homologie de 
$DHAB_*D^{odd(even)}({\Bbb Q})$. $\Box$

\vspace{2mm}

\noindent{\bf 2.4. Supercommutativit\'e des 
alg\`ebres (de Hopf) de l'homologie de}
{\bf $DHABD^{odd(even)}(\Bbbk )$ et de $DHAB_*D^{even}(\Bbbk )$}

\vspace{2mm}

{\bf Th\'eor\`eme 2.4.1.} {\it L'alg\`ebre (de Hopf) de l'homologie de
$DHABD^{odd(even)}(\Bbbk )$ est supercommutative pour n'importe quel 
anneau commutatif $\Bbbk$.} $\Box$

{\bf D\'emonstration du Th\'eor\`eme 2.4.1:} Le Th\'eor\`eme est une 
cons\'equence facile du Th\'eor\`eme 2.4.10. $\Box$

{\bf Th\'eor\`eme 2.4.2.} {\it L'alg\`ebre (de Hopf) de l'homologie
de $DHAB_*D^{even}(\Bbbk )$ est supercommutative pour n'importe
quel anneau commutatif $\Bbbk$.} $\Box$

{\bf D\'emonstration du Th\'eor\`eme 2.4.2:} Le Th\'eor\`eme est une 
cons\'equence facile du Th\'eor\`eme 2.4.40. $\Box$

Introduisons  des notations n\'ecessaires:

Soient $A$ et $B$ deux diagrammes de crochets (on cosid\`ere simultan\'ement
les cas de $d$ pair et impair).
Soient $t_1< t_2<\dots<
t_\beta$ les points du diagramme $B$.
 Alors on d\'efinit les diagrammes $A_{t_j}$, $1\leq
j\hm\leq\beta$, comme des diagrammes \'equivalents \`a $A$ et concentr\'es
dans de tr\`es petits voisinages des points $t_j$, $1\leq
j\hm\leq\beta$.

D\'efinissons
$$
A\vartriangleright_jB:\,=\begin{cases}
B|_{x_{t_j}=A_{t_j}}, & \text{si $1\le j\le\beta$;}\\
0, & \text{si $j>\beta$.}
\end{cases}
\eqno(2.4.3)$$
(La notation $B|_{x_{t_j}=A}$ est donn\'ee dans les sections 1.3.1, 1.3.2
pour $d$ impair et pair respectivement.)

D\'efinissons \'egalement
$$A\vartriangleright B:\,=\sum_{j=1}^\infty A\vartriangleright_jB.
\eqno(2.4.4)$$

Evidemment, si le bidegr\'e du diagramme $A$ est $(i_1,j_1)$, celui de
$B$ est $(i_2,j_2=\beta)$, alors le bidegr\'e des \'el\'ements
$A\vartriangleright_jB$, $A\vartriangleright B$ \'egale $(i_1+i_2,j_1+j_2-1)$.

On a imm\'ediatement le lemme suivant.

{\bf Lemme 2.4.5.} {\it Les op\'erations $\vartriangleright_j$,$j\in\Bbb N$,
$\vartriangleright$ introduites ci-dessus d\'efinissent correctement des 
applications
 $$
BD^{odd(even)}(\Bbbk)\otimes BD^{odd(even)}(\Bbbk)\rightarrow 
BD^{odd(even)}(\Bbbk).
\eqno(2.4.6)$$
($BD^{odd(even)}(\Bbbk )$ d\'esigne l'espace des $B$-diagrammes.) Ces 
applications sont de bidegr\'e $(0,-1)$.} $\Box$

{\bf Remarque 2.4.7.} Selon les Remarques 1.3.14, 1.3.35
$$
\partial A=(x_{t_-}\cdot x_{t+})\vartriangleright A-
(-1)^{p(A)-1}A\vartriangleright(x_{t_-}\cdot x_{t_+})
\eqno(2.4.8)$$
pour $d$ impair, o\`u $p(A)$ d\'esigne le poids de $A$ (voir la section 1.3.0),
$p(A)\equiv\tilde{A}\,{\rm mod}\,2$;
$$
\partial A=(x_{t-}\wedge x_{t+})\vartriangleright A-
(-1)^{p(A)-1}A\vartriangleright(x_{t_-}\wedge x_{t_+})
\eqno(2.4.9)$$
pour $d$ pair (dans ce cas $p(A)\equiv(\tilde{A}
+1)\,{\rm mod}\,2$). $\Box$

{\bf Th\'eor\`eme 2.4.10.} [T] {\it Pour tous deux diagrammes} $A,B\in
BD^{odd(even)}(\Bbbk)$
$$
\partial(A\vartriangleright B)=(\partial A)\vartriangleright B+(-1)^{p(A)-1}
A\vartriangleright(\partial B)+
(-1)^{p(A)-1}\left(A*B-(-1)^{p(A)p(B)}B*A\right).~\Box
\eqno(2.4.11)$$

Ici ``$*$'' signifie la multiplication dans $DHABD^{odd(even)}(\Bbbk )$
(voir la section 2.2).

{\bf Remarque 2.4.12.} L'affirmation du Th\'eor\`eme 2.4.10 est aussi
vraie pour les $B$-diagrammes g\'en\'eralis\'es. $\Box$

Maintenant consid\'erons $DHAB_*D^{odd(even)}(\Bbbk )$. La diff\'erentielle
$\partial$ sur $DHAB_*D^{odd(even)}(\Bbbk )$ peut \^etre d\'ecompos\'ee dans la
somme $\bar{\partial}+\,\bar{\!\bar\partial}$, o\`u $\bar{\partial}$
(resp. $\bar{\!\bar\partial}$) est la partie de $\partial$, qui garde
le nombre des ast\'erisques en augmentant en 1 le nombre
des composantes minimales (resp. diminue en 1 le nombre des ast\'erisques
tout en gardant le nombre des composantes minimales). Nous allons d\'esigner
 $DHAB_*D^{odd(even)}(\Bbbk )$, 
consid\'er\'ee comme une alg\`ebre de Hopf (on oublie la diff\'erentielle 
$\partial$),
par $HAB_*D^{odd(even)}(\Bbbk )$. Evidemment, l'op\'erateur 
$\bar{\partial}$ sers d'une diff\'erentielle pour cette alg\`ebre de Hopf.
L'alg\`ebre de Hopf diff\'erentielle ainsi obtenue sera d\'esign\'ee
par $(HAB_*D^{odd(even)}(\Bbbk );\bar{\partial})$.

{\bf Th\'eor\`eme 2.4.13.}[T] {\it L'alg\`ebre (de Hopf) de l'homologie
de  $(HAB_*D^{odd}(\Bbbk );\bar{\partial})$ 
(resp. $(HAB_*D^{even}(\Bbbk );\bar{\partial})$) 
est le produit tensoriel
de l'homologie de $DHABD^{odd}(\Bbbk )$ 
(resp. $DHABD^{even}(\Bbbk )$) 
sur l'alg\`ebre de Hopf
libre $F^{odd}(x)$ (resp. $F^{even}(x)$) \`a un g\'en\'erateur impair 
(resp. pair) primitif $x$, qui correspond au diagramme se composant 
d'un seul ast\'erisque isol\'e.} $\Box$

Notons que ce th\'eor\`eme entra\^\i ne la non-supercommutativit\'e
de l'alg\`ebre (de Hopf) de l'homologie de 
$(HAB_*D^{odd}(\Bbbk );\bar{\partial})$ (au contraire de 
$(HAB_*D^{even}(\Bbbk );\bar{\partial})$), parce que $F^{odd}(x)$ n'est 
pas supercommutative.

Soient $A$ et $B$ deux $*$-diagrammes de crochets; soient 
$t^{(*)}_1<t^{(*)}_2<\dots<t^{(*)}_\beta$ les points  du diagramme $B$ 
(qui peuvent avoir des 
ast\'erisques); 
$\tau^{(*)}_1<\tau^{(*)}_2<\dots <\tau^{(*)}_\alpha$ les points de $A$
(qui peuvent \'egalement avoir des 
ast\'erisques).
De mani\`ere analogue on d\'efinit les diagrammes $A_{t_j^{(*)}},$ $1\leq
j\leq\beta$, comme des diagrammes \'equivalents \`a $A$ et concentr\'es dans 
de tr\`es petits voisinages des points $t_j^{(*)}$.

On note par  $I_0$ (resp. $I_*$), 
$I_0\sqcup I_*=\{1,\dots,\alpha\}$, l'ensemble 
des index des points $\tau_i,$ $i\in I_0$ (resp. $\tau^*_i,$
$i\in I_*$), o\`u il n'y a pas d'ast\'erisques (resp. il y en a).

D\'efinissons
$$
A^*:=(-1)^{p(A)-1}\sum_{i\in I_0}A|_{x_{\tau_i}=x_{\tau_i^*}}.
\eqno(2.4.14)$$

Je vais rappeler, que $p(A)$ est le poids du diagramme $A$, 
$p(A)\equiv (\tilde{A} +1)\,{\rm mod}\,2$. Notons que $p(A^*)=p(A)+d-1\equiv
(p(A)+1)\,{\rm mod}\,2$.

{\bf Remarque 2.4.15.} $(A^*)^*=0$. $\Box$

D\'efinissons
$$
A\bar{\vartriangleright }_jB:=\begin{cases}B|_{x_{t_j}=A_{t_j}},& 
\text{si dans $t_j$,
$1\le j\le\beta$, il n'y a pas d'ast\'erisque;}\\
B|_{x_{t_j^*}=A^*_{t^*_j}},& \text{si dans $t^*_j$, $1\le j\le\beta$, 
il y a un ast\'erisque;}\\
0,& j>\beta .
\end{cases}
\eqno(2.4.16)$$

D\'efinissons \'egalement
$$
A\bar{\vartriangleright }B:=\sum_{j=1}^\infty A\bar{\vartriangleright }_jB.
\eqno(2.4.17)$$

{\bf Exemple 2.4.18.} Si dans tous les $\alpha$ points de $A$ et dans tous les
$\beta$ points de $B$ il y a des ast\'erisques, alors 
$A\bar{\vartriangleright }B=0.$ $\Box$

On a imm\'ediatement le lemme suivant.

{\bf Lemme 2.4.19.} {\it Les op\'erations $\bar{\vartriangleright }_j$,
$j\in\Bbb N$,
$\bar{\vartriangleright }$ d\'efinissent correctement des applications
$$
B_*D^{even}(\Bbbk)\otimes B_*D^{even}(\Bbbk)\rightarrow B_*D^{even}(\Bbbk).
\eqno(2.4.20)$$
($B_*D^{even}(\Bbbk )$ signifie l'espace des $*$-diagrammes de crochets pour 
$d$ pair) de bidegr\'e $(0,-1)$.} $\Box$

{\bf Remarque 2.4.21.} $
\bar{\partial} A=(x_{t-}\wedge x_{t+})\bar{\vartriangleright} A-
(-1)^{p(A)-1}A\bar{\vartriangleright}(x_{t_-}\wedge x_{t_+}).$ $\Box$

{\bf Th\'eor\`eme 2.4.22.}[T] {\it Pour tous deux diagrammes}
$A, B\in
B_*D^{even}(\Bbbk)$
$$
\bar{\partial}(A\bar{\vartriangleright }B)=(\bar{\partial}A)
\bar{\vartriangleright }B+(-1)^{p(A)-1}
A\bar{\vartriangleright }(\bar{\partial}B)+
(-1)^{p(A)-1}\left(A*B-(-1)^{p(A)p(B)}B*A\right).~\Box
\eqno(2.4.23)$$

Ici ``$*$'' signifie la multiplication dans $HAB_*D^{even}(\Bbbk)$.

{\bf Remarque 2.4.24.} Le Th\'eor\`eme 2.4.22 entra\^\i ne la 
supercommutativit\'e de l'alg\`ebre (de Hopf) de l'homologie de 
$(HAB_*D^{even}(\Bbbk),\bar{\partial}).~\Box$

Soit ${\gothg}$ une superalg\`ebre de Lie. D\'efinissons un op\'erateur
$$
\delta:\Lambda^*{\gothg}\rightarrow\Lambda^{*-1}{\gothg},
\eqno(2.4.25)$$
en prenant pour $A=A_1\wedge \dots\wedge A_k,$ $A_i\in{\gothg},$
$$
\delta(A):=\delta(A_1\wedge\dots\wedge A_k)=\sum_{i<j}(-
1)^{\lambda_{ij}+\tilde A_i}[A_i,A_j]\wedge
A_1\wedge\dots\hat{A}_i\dots\hat{A}_j\dots\wedge A_k,
\eqno(2.4.26)$$
o\`u $\lambda_{ij}=(\tilde A_i+1)\sum\limits_{p=1}^{i-1}(\tilde
A_p+1)+(\tilde A_j+1)\sum\limits_{q=1\atop q\neq i}^{j-1}(\tilde A_q+1)$.

Il est facile de voir, que $\delta$ peut \^etre correctement
prolong\'e sur toute l'alg\`ebre ext\'erieure $\Lambda^*{\gothg}$.
Notons \'egalement que $\delta^2=0.$ La couple 
$(\Lambda^*{\gothg},\delta)$ d\'efinit le dit {\it  complexe
de cha\^\i nes, associ\'e \`a ${\gothg}$} et calculant son homologie.

Pour tous  \'el\'ements purs $A,B\in\Lambda^*{\gothg}$ on a 
les formules suivantes
$$
\delta(A\wedge B)=\delta(A)\wedge B+(-1)^{\tilde A+1}A\wedge
\delta(B)+(-1)^{\tilde A}[A,B];
\eqno(2.4.27)$$
$$
\delta([A,B])=[\delta(A),B]+(-1)^{\tilde A}[A,\delta(B)].
\eqno(2.4.28)$$
O\`u $[.,.]$ est le crochet de Schouten (voir la section 1.3.2).

De mani\`ere analogue (selon la formule (2.4.27), o\`u les $A_i$ sont
des crochets correspondant aux composantes minimales) on d\'efinit 
l'op\'erateur
$$
\delta:B_*D^{even}(\Bbbk)\rightarrow B_*D^{even}(\Bbbk)
\eqno(2.4.29)$$
sur l'espace $B_*D^{even}(\Bbbk)$. $\delta$ est de bidegr\'e $(1,0)$ et 
donc impair.

Pour le diagramme $A$ on d\'efinit
$$
A^0:=(-1)^{p(A)-1}\delta(A)=(-1)^{\tilde A}\delta(A);
\eqno(2.4.30)$$
$$
A^\star:=A^*+A^0.
\eqno(2.4.31)$$

{\bf Remarque 2.4.32.} $(A^\star)^\star=0.$ $\Box$

D\'efinissons

$$
A\bar{\bar{\vartriangleright }}_jB:=\begin{cases}
B|_{x_{t^*_j}=A^0_{t^*_j}}, & \text{si $1\le j\le\beta$ et dans $t_j^*$ il y
a un ast\'erisque;}\\
0,& \text{si non;}
\end{cases}
\eqno(2.4.33)$$
$$
A\vartriangleright _jB:=A\bar{\vartriangleright }_jB+
A\bar{\bar{\vartriangleright }}_jB;
\eqno(2.4.34)$$
$$
A\bar{\bar{\vartriangleright }}B:=
\sum_{j=1}^\infty A\bar{\bar{\vartriangleright }}_jB;
\eqno(2.4.35)$$
$$
A\vartriangleright B:=A\bar{\vartriangleright }B+
A\bar{\bar{\vartriangleright }}B.
\eqno(2.4.36)$$

{\bf Lemme 2.4.37.} {\it Les op\'erations $\bar{\bar{\vartriangleright}}_j,$
$\vartriangleright_j,$ $\bar{\bar{\vartriangleright}},$
$\vartriangleright$ d\'efinissent correctement des applications
$$
B_*D^{even}(\Bbbk)\otimes B_*D^{even}(\Bbbk)\rightarrow B_*D^{even}(\Bbbk)
\eqno(2.4.38)$$
de bidegr\'e $(0,-1)$.} $\Box$

{\bf D\'emonstration du Lemme 2.4.37:} \'evident. $\Box$

{\bf Remarque 2.4.39.} $
\partial A=(x_{t-}\wedge x_{t+})\vartriangleright A-
(-1)^{p(A)-1}A\vartriangleright(x_{t_-}\wedge x_{t_+}).$ $\Box$

{\bf Th\'eor\`eme 2.4.40.}[T] {\it Pour tous deux diagrammes} 
$A,B\in B_*D^{even}(\Bbbk)$
$$
\partial(A\vartriangleright B)=(\partial A)\vartriangleright B+(-1)^{p(A)-1}
A\vartriangleright(\partial B)+
(-1)^{p(A)-1}\left(A*B-(-1)^{p(A)p(B)}B*A\right).~\Box
\eqno(2.4.41)$$

{\bf Remarque 2.4.42.} Pour le cas o\`u $d$ est impair l'alg\`ebre 
(de Hopf) de
l'homologie de $(HAB_*D^{odd}(\Bbbk );\bar{\partial})$ n'est pas
supercommutative, 
alors on ne peut pas \'ecrire
m\^eme le terme grossier (analogue \`a (2.4.18)) de l'homotopie
du supercommutateur dans $DHAB_*D^{odd}(\Bbbk )$. $\Box$

\vspace{4mm}

\noindent{\bf 2.5. L'homologie de  $DHAB_*D^{odd(even)}(\Bbbk )$ et de
$DHAB_*D^{even}(\Bbbk )$ 
forme des alg\`ebres de Gerstenhaber}

\vspace{2mm}

Il est facile de montrer, que les complexes $CBD^{odd(even)}(\Bbbk )$,
$CB_*D^{even}(\Bbbk )$ munis du crochet impair $\{.,.\}$:
$$
\{A,B\}:=A\vartriangleright B-(-1)^{(p(A)-1)(p(B)-1)}B
\vartriangleright A
\eqno(2.5.1)
$$
sont en fait des superalg\`ebres de Lie diff\'erentielles, o\`u la 
diff\'erentielle $\partial$ (selon les Remarques 2.4.7, 2.4.39) est le 
crochet avec $x_{t_-}\cdot x_{t_+}$ (ou $x_{t_-}\wedge x_{t_+}$):
$$
\partial A=\{ x_{t_-}\cdot x_{t_+},A\}\, {\rm ou}\,\{ x_{t_-}
\wedge x_{t_+},A\}.$$

Dans l'homologie le crochet $\{.,.\}$ sera concordant avec la 
muliplication (qui y est supercommutative), c'est-\`a-dire pour toutes trois
classes $x$, $y$, $z$ de l'homologie on a
$$\{ x,yz\} =\{x,y\}z+(-1)^{(p(x)-1)p(y)}y\{x,z\}.
\eqno(2.5.2)
$$

Ce fait suit des consid\'erations g\'en\'erales pour les op\'erades,
voir la section 3 et le Th\'eor\`eme 3.2.4 en particulier.

{\bf D\'efinition 2.5.3.}
Les alg\`ebres supercommutatives munis d'un crochet de Lie impair, qui est
compatible avec la multiplication selon (2.5.2), sont appel\'ees
{\it alg\`ebres de Gerstenhaber} (dans cette situation ce crochet-l\`a est 
appel\'e {\it crochet (impair) de Kirillov}). $\Box$

{\bf Exemple 2.5.4.} L'alg\`ebre ext\'erieure d'une superalg\`ebre de Lie
avec le crochet de Schouten (voir la section 1.3.2) fait une alg\`ebre de 
Gerstenhaber. $\Box$

Notons que le crochet $\{.,.\}$ conserve l'espace des \'el\'ements primitifs.
Alors, si $\Bbbk$ est un corps de caract\'eristique nulle, l'alg\`ebre de
Gerstenhaber inconnue est donc l'alg\`ebre ext\'erieure de la superalg\`ebre 
de Lie (avec le crochet $\{.,.\}$) des \'el\'ements primitifs de l'homologie.

\vspace{4mm}

\noindent{\bf 2.6. Calcul de l'homologie de $DHABD^{odd(even)}({\Bbb Z})$
pour la complexit\'e $i\le 2$}

\vspace{2mm}

\noindent{\bf 2.6.1. Cas o\`u $d$ est pair}

{\bf Notation 2.6.1.} Dans cette section et dans celle qui la suivent le
g\'en\'erateur $x_i$ correspond toujours au $i$-i\`eme point, compt\'e
de gauche \`a droite dans les diagrammes correspondants. $\Box$

$\underline{i=0}$. On n'a que le diagramme trivial; ce diagramme ne
vient pas de la g\'eom\'etrie du discriminant, mais on l'ajoute dans tous
les complexes $CB_{(*)}D^{odd(even)}(\Bbbk )$, parce qu'il correspond \`a 
l'unit\'e (perdue sous la dualit\'e d'Alexander) de l'homologie.
Ce diagramme d\'efinit un cycle \'etant l'unit\'e dans l'alg\`ebre (de Hopf)
de l'homologie de $DHABD^{odd(even)}({\Bbbk})$.

$\underline{i=1}$. On a aussi un seul $B$-diagramme --- $[x_1x_2]$,
qui d\'efinit un g\'en\'erateur
 impair (primitif) $y$ de bidegr\'e $(1,2)$ dans l'homologie de 
$DHABD^{even}({\Bbb Z})$.

$\underline{i=2}$. L'espace des diagrammes de bidegr\'e (2,3) est de dimension
2. Pour une base on prend
$$\begin{array}{l}
a_1=[[x_1x_2]x_3];\\
a_2=[[x_1x_3]x_2].
\end{array}
\eqno(2.6.2)$$

L'espace des diagrammes de bidegr\'e (2,4)
est de dimension 3. Pour une base on prend
$$
\begin{array}{l}
b_1=[x_1x_4]\wedge[x_2x_3];\\
b_2=[x_1x_3]\wedge[x_2x_4];\\
b_3=[x_1x_2]\wedge[x_3x_4].
\end{array}
\eqno(2.6.3)$$

Trouvons $\partial a_1,$ $\partial a_2$:
$$
\partial a_1=\partial([[x_1x_2]x_3])=P
\bigl( [[x_1\wedge x_2,x_3]x_4]+[[x_1,x_2\wedge
x_3]x_4]\bigr) = $$ $$
P\bigl( [x_1\wedge[x_2x_3],x_4]-
[x_2\wedge[x_1x_3],x_4]+[[x_1x_2]\wedge x_3,x_4]-
[[x_1x_3]\wedge x_2,x_4]\bigr)=
$$
$$[x_1x_4]\wedge[x_2x_3]+[x_1x_3]\wedge[x_2x_4]+
[x_1x_2]\wedge[x_3x_4]-[x_1x_3]\wedge[x_2x_4]=b_1+b_3.
\eqno(2.6.4)$$
$$
\partial a_2=\partial([[x_1x_3]x_2])=P\bigl( [[x_1\wedge x_2,x_4]x_3]+
[[x_1,x_3\wedge x_4],x_2]\bigr) =
$$
$$P\bigl(
[x_1\wedge[x_2x_4],x_3]-
[x_2\wedge[x_1x_4],x_3]+[[x_1x_3]\wedge x_4,x_2]-[[x_1x_4]\wedge x_3,x_2]
\bigr) =
$$
$$
[x_1x_3]\wedge[x_2x_4]+[x_1x_4]\wedge[x_2x_3]-[x_1x_3]\wedge[x_2x_4]+
[x_1x_4]\wedge[x_2x_3]=2b_1.
\eqno(2.6.5)$$

D'o\`u on d\'eduit, que dans le bidegr\'e $(2,3)$ l'homologie est triviale;
dans le bidegr\'e $(2,4)$ l'homologie est isomorphe \`a 
${\Bbb Z}\oplus{\Bbb Z}_2$. La torsion ${\Bbb Z}_2$ est 
engendr\'e par le cycle $y^2=y*y=b_3=[x_1x_2]\wedge[x_3x_4].$
Pour un $\Bbb Z$-g\'en\'erateur on peut prendre $u=b_2=[x_1x_3]\wedge[x_2x_4].$
Il est facile de voir, que $y^2$ et $u$ sont primitifs.

\vspace{2mm}

\noindent{\bf 2.6.2. Cas o\`u $d$ est impair}

$\underline{i=0}$. On n'a que le diagramme trivial d\'efinissant 
l'unit\'e dans l'alg\`ebre (de Hopf) de l'homologie.

$\underline{i=1}$. On n'a que le diagramme $[x_1x_2]$, qui d\'efinit un 
g\'en\'erateur pair (primitif) $y$ de bidegr\'e $(1,2)$ dans l'homologie de 
$DHABD^{odd}(\Bbb Z)$.

$\underline{i=2}.$ Pour une base dans l'espace des diagrammes on prend 
$a_1,a_2;b_1,b_2,b_3$ respectivement d\'efinis par les formules
(2.6.1), (2.6.2) (on remplace la multiplication ext\'erieure par la 
multiplication normale).

Trouvons $\partial a_1,\partial a_2$:
$$
\partial a_1=P\bigl( [[x_1\cdot x_2,x_3],x_4]-[[x_1,x_2\cdot x_3],x_4]
\bigr)=
$$
$$
P\bigl([x_1\cdot[x_2x_3],x_4]-[x_2\cdot[x_1x_3],x_4]-[[x_1x_2]\cdot x_3,x_4]+
[[x_1x_3]\cdot x_2,x_4]\bigr)=
$$
$$
[x_1x_4]\cdot[x_2x_3]-[x_1x_3]\cdot[x_2x_4]-
[x_1x_2]\cdot[x_3x_4]+[x_1x_3]\cdot[x_2x_4]=b_1-b_3.
\eqno(2.6.6)$$
$$
\partial a_2=P\bigl([[x_1\cdot x_2,x_4],x_3]-[[x_1,x_3\cdot x_4],x_2]\bigr)=
$$
$$
P\bigl([x_1\cdot[x_2x_4],x_3]-[x_2\cdot[x_1x_4],x_3]-[[x_1x_3]\cdot x_4,x_2]+
[[x_1x_4]\cdot x_3,x_2]\bigr)=
$$
$$
[x_1x_3]\cdot[x_2x_4]-[x_1x_4]\cdot[x_2x_3]-
[x_1x_3]\cdot[x_2x_4]+[x_1x_4]\cdot[x_2x_3]=0.
\eqno(2.6.7)$$

D'o\`u on d\'eduit, que dans le bidegr\'e (2,3) l'homologie
est isomorphe \`a $\Bbb Z$; il appara\^\i t un nouveau g\'en\'erateur
impair (primitif) $z=[[x_1x_3]x_2]$. Il est facile de voir, que
$\{y,y\}=\{[x_1x_2],[x_1x_2]\}=-2z$. Dans le bidegr\'e $(2,4)$ 
l'homologie est isomorphe \`a ${\Bbb Z}\oplus{\Bbb Z}$, pour une base
on peut prendre $y^2=b_3=[x_1x_2]\cdot[x_3x_4]$ et 
$u=b_2=[x_1x_3]\cdot[x_2x_4].$ Au lieu de $u$ on peut prendre l'\'el\'ement
$$
u'=\log_\star id(u)=u-y^2,
\eqno(2.6.8)$$
qui est bien primitif (la d\'escription de l'op\'erateur $\log_\star id$
est donn\'ee dans la section 2.1).

\vspace{4mm}

\noindent{\bf 2.7. Rapport entre les alg\`ebres de Hopf de l'homologie
de
$DHABD^{odd(even)}(\Bbb Q)$
et de  $DHAB_*D^{odd(even)}(\Bbb Q)$ }

\vspace{2mm}

\noindent{\bf 2.7.1. Cas o\`u $d$ est pair}

Consid\'erons le diagramme commutatif (2.2.7) dans le cas o\`u $d$ est pair, 
pour ${\Bbbk}={\Bbb Q}$ (ou bien pour n'importe quel
autre corps de caract\'eristique nulle).

{\bf Th\'eor\`eme 2.7.1.}[T] {\it Si $d$ est pair et $\Bbbk$ est
 un corps de caract\'eristique nulle, 
l'application $i_1$ (resp. $p_1$) du diagramme commutatif 
(2.2.7) induit une application surjective en homologie, dont le noyau
est l'id\'eal engendr\'e par $y=[x_1,x_2]$ (g\'en\'erateur pair primitif --
voir la section 2.6).} $\Box$

Notons que l'application $i_1$ en homologie est compl\`etement
d\'efinie par sa restriction sur l'espace des \'el\'ements primitifs. 
Les espaces des \'el\'ements primitifs de 
 $DHABD^{even}({\Bbb Q})$ et de $DHAB_*D^{even}({\Bbb Q})$ poss\`edent
le crochet $\{.,.\}$ (voir la section 2.5). L'application $i_1$
d\'efinit un morphisme des 
alg\`ebres de Lie en question. On en conclut imm\'ediatement, que $y=[x_1,x_2]$
est un \'el\'ement central de l'alg\`ebre de Lie
(avec le crochet  $\{.,.\}$) des \'el\'ements primitifs de
l'homologie de $DHABD^{even}({\Bbb Q})$.

\vspace{2mm}

\noindent{\bf 2.7.2. Cas o\`u  $d$ est impair}

{\bf Th\'eor\`eme 2.7.2.}[T] {\it Si $d$ est impair
et $\Bbbk$ est un corps de caract\'eristique nulle,
l'application $i_1$ (resp. $p_1$) du diagramme commutatif 
(2.2.7) induit une application surjective en homologie, dont le noyau
est l'id\'eal engendr\'e par les g\'en\'erateurs
primitifs $y=[x_1,x_2]$ (qui est pair) et $z=[[x_1x_3]x_2]$
(qui est impair), voir la section 2.6.} $\Box$

Ce th\'eor\`eme entra\^\i ne la supercommutativit\'e de l'homologie de 
$DHAB_*D^{odd}({\Bbb Q})$

\vspace{4mm}

\noindent{\bf 2.8. Big\`ebres des (super)diagrammes de cordes}

\vspace{2mm}

\noindent{\bf 2.8.1. Cas o\`u $d$ est impair. 
Big\`ebres des diagrammes de cordes}

On va d\'esigner par  $H_{*,*}(DHAB_{(0)}D^{odd}(\Bbbk))$
l'alg\`ebre (de Hopf) bigradu\'ee de l'homologie de $DHAB_{(0)}D^{odd}(\Bbbk)$.
Il est facile de montrer, que (voir les Exemples 1.3.17, 1.3.20)
$$
BChD^{odd}(\Bbbk):=\oplus_{i=0}^{+\infty}H_{i,2i}(DHABD^{odd}(\Bbbk));
\eqno(2.8.1)$$
$$
BCh_0D^{odd}(\Bbbk):=\oplus_{i=0}^{+\infty}H_{i,2i}(DHAB_0D^{odd}(\Bbbk))
\eqno(2.8.2)$$
sont des big\`ebres pour n'importe quel anneau commutatif $\Bbbk$.
Ces big\`ebres sont les dites {\it big\`ebres des diagrammes de cordes}
($BChD$ vient de l'anglais ``Bialgebra of Chord Diagrams'').

L'application naturelle
$$
p:BChD^{odd}(\Bbbk)\rightarrow BCh_0D^{odd}(\Bbbk)
\eqno(2.8.3)$$
est surjective, dont le noyau est l'id\'eal, engendr\'e par le g\'en\'erateur
primitif correspondant \`a une corde isol\'ee $[x_1,x_2]$. Ces big\`ebres sont
commutatives, parce que le bidegr\'e $(i,2i)$ est toujours pair.

Normalement on consid\`ere les big\`ebres des diagrammes (de cordes)
sur un cercle $S^1$, c'est parce que les th\'eories des invariants (de type 
fini -- en particulier) de n\oe uds et de n\oe uds non-compacts
co\"\i ncident. Il est facile de montrer, que tout diagramme de cordes $D$
avec les points $t_1<t_2<\dots<t_{2n}$ sur ${\Bbb R}$ est \'egal
(grace aux relations de quatre termes) au diagramme, qui est obtenu de $D$
par la permutation circulaire des points:
$$
(t'_1,t'_2,\dots,t'_{2n})=(t_2,\dots,t_{2n},t_1).
\eqno(2.8.4)$$

La big\`ebre des diagrammes de cordes a \'et\'e intens\'ement \'etudi\'ee
pendant les derni\`eres ann\'ees, voir [BN, ChD, ChDL, K1, Kn, L, NS, S, Z]. 
Il est connu que jusqu'\`a la complexit\'e $i=12$ il n'y a pas de torsion dans
$BChD(\Bbb Z)$ (voir [K]), ce qui confirme indirectement la conjecture
de Vassiliev sur la trivialit\'e homotopique de la filtration (0.1.2) 
dans le discriminant r\'esolu $\sigma$. Comme S.K.Lando l'a d\'emontr\'e,
voir [L], $BChD({\Bbb Z})$ est engendr\'ee par ses \'el\'ements primitifs. 
Alors la big\`ebre $BChD({\Bbb Z})$, \'etant quotient\'ee par les torsions 
possibles, est l'alg\`ebre sym\'etrique sur $\Bbb Z$ de l'espace de ses 
\'el\'ements primitifs (comme elle est commutative aussi bien que 
cocommutative). Les caculs, commenc\'es par V.A.Vassiliev
et poursuits par D.Bar-Natan [BN] et
par J.A.Kneissler [Kn], donnent les valeurs suivantes des dimensions
$p_i$, $i=1,\dots,12,$ des espaces des \'el\'ements primitifs de complexit\'e
$i$.

\vspace{5mm}

\quad\quad
\quad\quad\quad\begin{tabular}{c|c|c|c|c|c|c|c|c|c|c|c|c}
$i$ & 1 & 2 & 3 & 4 & 5 & 6 & 7 & 8 & 9 & 10 & 11 & 12\\
\hline
$p_i$ & 1 & 1 & 1 & 2 & 3 & 5 & 8 & 12 & 18 & 27 & 39 & 55\\
\end{tabular}
\quad\quad\quad\quad
\quad\quad\quad\quad
(2.8.5)

\vspace{5mm}

\noindent{\bf 2.8.2. Cas o\`u $d$ est pair. 
Big\`ebres des superdiagrammes de cordes}

De mani\`eres analogue on d\'efinit les {\it big\`ebres des superdiagrammes 
de cordes} (voir l'Exemple 1.3.36):
$$
BChD^{even}(\Bbbk):=\oplus_{i=0}^{+\infty}H_{i,2i}(DHABD^{even}(\Bbbk));
\eqno(2.8.6)$$
$$
BCh_0D^{even}(\Bbbk):=\oplus_{i=0}^{+\infty}H_{i,2i}(DHAB_0D^{even}(\Bbbk)).
\eqno(2.8.7)$$

L'application naturelle 
$$
p:BChD^{even}(\Bbbk)\rightarrow BCh_0D^{even}(\Bbbk)
\eqno(2.8.8)$$
est \'egalement surjective, dont le noyau est l'id\'eal engendr\'e
par le g\'en\'erateur impair primitif correspondant \`a une corde isol\'ee
$[x_1,x_2].$

Dans ce cas-l\`a le diagramme $D$ avec les points $t_1<\dots<t_{2n}$ 
sur ${\Bbb R}$ n'est plus \'equivalent au diagramme $D'$ obtenu de $D$
par la permutation circulaire des ponts (2.8.4)
(par exemple, c'est pas vrai pour une corde isol\'ee $[x_1x_2]$),
mais par contre tous deux diagrammes semblables sont toujours 
\'egaux dans $BCh_0D^{even}(\Bbbk)$. Alors $BCh_0D^{even}(\Bbbk)$ peut \^etre
vue comme l'espace des superdiagrammes de cordes sur un cercle $S^1$
quotient\'e par les relations de quatre termes et par les relations d'un terme.

L'espace des superdiagrammes de cordes sur un cercle, quotient\'e
seulement par les relations de quatre termes est isomorphe 
\`a l'espace des superdiagrammes de cordes sur une droite  quotient\'e
par les relations de quatre termes et par les \underline{doubles} relations 
d'un terme.

\vspace{4mm}

\noindent{\bf 2.9. Calcul de l'homologie de $DHAB_0D^{odd(even)}(\Bbb Z)$ 
pour la
complexit\'e $i\le 3$}

\vspace{2mm}

\noindent{\bf 2.9.1. Cas o\`u $d$ est pair}

$\underline{i=0}.$ On n'a que le diagramme trivial, qui d\'efinit 
l'unit\'e dans l'alg\`ebre (de Hopf) de l'homologie.

$\underline{i=1}.$ L'espace des $B_0$-diagrammes y est trivial.

$\underline{i=2}.$ On n'a qu'un seul diagramme non-nul
$[x_1x_3]\wedge[x_2x_4]$, qui d\'efinit un g\'en\'erateur pair (primitif)
de bidegr\'e $(2,4)$ dans l'alg\`ebre (de Hopf) de l'homologie de 
$DHAB_0D^{even}(\Bbb Z)$.

$\underline{i=3}.$ Si la deuxi\`eme graduation $j=4$, alors de tels 
diagrammes ont une composante minimale. L'espace de ces diagrammes est 
isomorphe \`a $\Bbb Z$. Pour un g\'en\'erateur on prendra 
$$
a_1=[[x_1x_3],[x_2x_4]].
\eqno(2.9.1)$$

Si $j=5$ pour une base dans l'espace des diagrammes on prend
$$
\left.\begin{array}{l}
b_1=[[x_2x_4]x_5]\wedge[x_1x_3];\\
b_2=[[x_2x_5]x_3]\wedge[x_1x_4];\\
b_3=[[x_1x_3]x_5]\wedge[x_2x_4];\\
b_4=[[x_1x_5]x_3]\wedge[x_2x_4];\\
b_5=[[x_1x_3]x_4]\wedge[x_2x_5];\\
b_6=[[x_1x_4]x_2]\wedge[x_3x_5].
\end{array}\right\}
\eqno(2.9.2)$$

Si $j=6$ pour une base on prend
$$
\left.\begin{array}{l}
c_1=[x_1x_6]\wedge[x_2x_4]\wedge[x_3x_5]\\
c_2=[x_1x_5]\wedge[x_2x_4]\wedge[x_3x_6]\\
c_3=[x_1x_4]\wedge[x_2x_6]\wedge[x_3x_5]\\
c_4=[x_1x_4]\wedge[x_2x_5]\wedge[x_3x_6]\\
c_5=[x_1x_3]\wedge[x_2x_5]\wedge[x_4x_6]
\end{array}\right\}
\eqno(2.9.3)$$

Trouvons
$$
\partial a_1=\partial([[x_1x_3],[x_2x_4]])=
$$
$$
=P([[x_1\wedge x_2,x_4],[x_3x_5]]+[[x_1,x_3\wedge x_4],[x_2x_5]]+
$$
$$
+[[x_1x_4],[x_2\wedge x_3,x_5]]+[[x_1x_3],[x_2,x_4\wedge x_5]])=
$$
$$
=[[x_2x_4]x_5]\wedge[x_1x_3]-[[x_1x_5]x_3]\wedge [x_2x_4]+
$$
$$
+[[x_1x_4]x_2]\wedge[x_3x_5]=b_1-b_4+b_6.
\eqno(2.9.4)$$

Trouvons
$$
\partial b_1=\partial([[x_2x_4]x_5]\wedge[x_1x_3])=
$$
$$
=P\left([[x_2\wedge x_3,x_5],x_6]\wedge[x_1x_4]+[[x_2,x_4\wedge
x_5],x_6]\wedge[x_1x_3]\right)=
$$
$$
=[x_1x_4]\wedge[x_2x_6]\wedge[x_3x_5]+[x_1x_4]\wedge[x_2x_5]\wedge[x_3x_6]+
$$
$$
+[x_1x_3]\wedge[x_2x_4]\wedge[x_5x_6]-[x_1x_3]\wedge[x_2x_5]\wedge[x_4x_6]=
$$
$$
=c_3+c_4-c_5.
\eqno(2.9.5)$$

De mani\`ere analogue on obtient
$$
\left.
\begin{array}{l}
\partial b_2=c_2+c_3-c_4\\
\partial b_3=c_1+c_3+c_4-c_5\\
\partial b_4=c_2+c_3\\
\partial b_5=c_2+c_3+c_4\\
\partial b_6=c_2-c_4+c_5
\end{array}\right\}
\eqno(2.9.6)$$

Alors les matrices de la diff\'erentielle $\partial$ seront suivantes:
$$
\left(\begin{array}{cccccc}
0 & 0 & 1 & 0 & 0 & 0\\
0 & 1 & 0 & 1 & 1 & 1\\
1 & 1 & 1 & 1 & 1 & 0\\
1 &-1 & 1 & 0 & 1 & -1 \\
-1 & 0 & -1& 0 & 0 & 1
\end{array}\right)
\left(\begin{array}{c}
1\\
0\\
0\\
-1\\
0\\
1
\end{array}\right).
\eqno(2.9.7)$$

D'o\`u on a:

Dans le bidegr\'e $(3,4)$ l'homologie est triviale.

Dans le bidegr\'e $(3,5)$ l'homologie est isomorphe \`a $\Bbb Z$. On obtient un 
g\'en\'erateur pair (primitif)
$$
w=b_2-2b_4+b_5=[[x_2x_5]x_3]\wedge[x_1x_4]-
$$
$$
-2[[x_1x_5]x_3]\wedge[x_2x_4]+[[x_1x_3]x_4]\wedge[x_2x_5].
\eqno(2.9.8)$$

{\bf Remarque 2.9.9.} Pour l'\'el\'ement
$$
w'=[[x_2x_5]x_3]\wedge[x_1x_4]-2[[x_1x_5]x_3]\wedge[x_2x_4]+
[[x_1x_4]x_3]\wedge[x_2x_5]
\eqno(2.9.10)$$
$\partial w'=0$ dans $DHABD^{even}(\Bbb Z)$. D'un autre c\^ot\'e
$w'=w$ dans $DHAB_0D^{even}(\Bbb Z)$, alors l'\'el\'ement $w$
(consid\'er\'e comme une classe de l'homologie de $DHAB_0D^{even}(\Bbb Z)$)
appartient \`a ${\rm Im}(p_1)_*$, o\`u $(p_1)_*$
est l'application en homologie induite par la projection $p_1$
(voir le diagramme 2.2.7). Il est facile de voir, que $w'$
(comme une classe de l'homologie) ne peut pas \^etre obtenu des 
diagrammes de cordes par des op\'erations de la multiplication $*$
et du crochet $\{.,.\}$. $\Box$

Dans le bidegr\'e $(3,6)$ on a les relations suivantes entre les 
$c_i$,
$i=1,\dots,5$.
$$
\left.\begin{array}{l}
c_1=c_4=0;\\
c_2=-c_3=-c_5.
\end{array}\right\}
\eqno(2.9.11)$$

Cela donne encore un g\'en\'erateur impair (primitif).

\vspace{2mm}

\noindent{\bf 2.9.2. Cas o\`u $d$ est impair}

$\underline{i=0}.$ On n'a que le diagramme trivial d\'efinissant
l'unit\'e dans l'alg\`ebre (de Hopf) de l'homologie.

$\underline{i=1}.$ L'espace des $B_0$-diagrammes est trivial.

$\underline{i=2}$. On n'a qu'un seul diagramme non-nul
$[x_1x_3]\cdot[x_2x_4]$, qui d\'efinit un g\'en\'erateur pair (primitif) de 
bidegr\'e $(2,4)$ dans l'alg\`ebre (de Hopf) de l'homologie de 
$DHAB_0D^{odd}(\Bbb Z).$

$\underline{i=3}.$ Comme une base pour $j=4,5,6$ on prend
respectivement les m\^emes $a_1;b_1,\dots,b_6;c_1,\dots,c_5$ en 
rempla\c{c}ant la multiplication ext\'erieure par la multiplication
normale.

Trouvons
$$
\partial a_1=\partial([[x_1x_3],[x_2x_4]])=
$$
$$
=P([[x_1\cdot x_2,x_4],[x_3x_5]]-[[x_1,x_3\cdot x_4],[x_2x_5]]+
$$
$$
+[[x_1x_4],[x_2\cdot x_3,x_5]]-[[x_1x_3],[x_2,x_4\cdot x_5]])=
$$
$$
=-[[x_2x_4]x_5]\cdot [x_1x_3]-2[[x_2x_5]x_3]\cdot [x_1x_4]+
$$
$$
+[[x_1x_5]x_3]\cdot [x_2x_4]+2[[x_1x_3]x_4]\cdot [x_2x_5]+
$$
$$
+[[x_1x_4]x_2]\cdot [x_3x_5]=-b_1-2b_2+b_4+2b_5+b_6.
\eqno(2.9.12)$$

Trouvons
$$
\partial b_1=\partial([[x_2x_4]x_5]\cdot [x_1x_3])=
$$
$$
=P\left([[x_2\cdot x_3,x_5],x_6]\cdot [x_1x_4]-[[x_2,x_4\cdot x_5],x_6]\cdot
[x_1x_3]\right)=
$$
$$
=[x_1x_4]\cdot [x_2x_6]\cdot [x_3x_5]-[x_1x_4]\cdot [x_2x_5]\cdot [x_3x_6]+
$$
$$
+[x_1x_3]\cdot [x_2x_4]\cdot [x_5x_6]+[x_1x_3]\cdot [x_2x_5]\cdot [x_4x_6]=
c_3-c_4+c_5.
\eqno(2.9.13)$$

De mani\`ere analogue on obtient
$$
\left.\begin{array}{l}
\partial b_2=c_2+c_3-c_4;\\
\partial b_3=c_1-c_3+c_4-c_5;\\
\partial b_4=-c_2+c_3;\\
\partial b_5=c_2+c_3-c_4;\\
\partial b_6=c_2-c_4+c_5.
\end{array}\right\}
\eqno(2.9.14)$$

Alors les matrices de la diff\'erentielle $\partial$ sont les suivantes:

$$
\left(\begin{array}{cccccc}
0 & 0 & 1 & 0 & 0 & 0\\
0 & 1 & 0 & -1 & 1 & 1\\
1 & 1 & -1 & 1 & 1 & 0\\
-1 & -1 & 1 & 0 & -1 & -1\\
1 & 0 & -1& 0 & 0 & 1
\end{array}\right)
\left(\begin{array}{c}
-1\\
-2\\
0\\
1\\
2\\
1
\end{array}\right).
\eqno(2.9.15)$$

Donc:

Dans le bidegr\'e $(3,4)$ l'homologie est triviale.

Dans le bidegr\'e $(3,5)$ l'homologie est isomorphe \`a $\Bbb Z$.
On obtient un g\'en\'erateur impair (primitif)
$$
w=b_2-b_5=[[x_2x_5]x_3]\cdot[x_1x_4]-[[x_1x_3]x_4]\cdot[x_2x_5].
\eqno(2.9.16)$$

{\bf Remarque 2.9.17.} Pour l'\'el\'ement
$$
w'=[[x_2x_5]x_3]\cdot[x_1x_4]+[[x_1x_4]x_3]\cdot[x_2x_5]
\eqno(2.9.18)$$
$\partial w'=0$ dans $DHABD^{odd}({\Bbb Z})$.
D'un autre c\^ot\'e, $w'=w$ dans $DHAB_0D^{odd}({\Bbb Z})$, alors
l'\'el\'ement $w$ (consid\'er\'e comme une classe de l'homologie
de $DHAB_0D^{odd}({\Bbb
Z})$) appartient \`a ${\rm Im}(p_1)_*$, o\`u $(p_1)_*$ est l'application
en homologie induite par la projection $p_1$ (voir le diagramme (2.2.7)).
Notons que
$$
\{ [x_1x_2],[x_1x_3]\cdot [x_2x_4]\}=[x_1x_2]\vartriangleright
([x_1x_3]\cdot [x_2x_4])+([x_1x_3]\cdot [x_2x_4])
\vartriangleright[x_1x_2]=$$
$$
=-[x_1x_3]\cdot[[x_2x_5]x_4]-[x_2x_4][[x_1x_5]x_3]-
[x_3x_5]\cdot [[x_1x_4]x_2]-[x_2x_5]\cdot [x_1[x_3x_4]]+
[x_1x_4]\cdot [[x_2x_3]x_5].
\eqno(2.9.19)$$
Dans $DHAB_0D^{odd}({\Bbb Z})$ cet \'el\'ement est \'egal \`a
$$
[x_1x_3]\cdot [[x_2x_4]x_5]-[x_2x_4]\cdot [[x_1x_5]x_3]-
[x_3x_5]\cdot [[x_1x_4]x_2]=$$
$$
=b_1-b_4-b_6=-\partial a_1-2b_2 +2b_5=-\partial a_1 -2w.
\eqno(2.9.20)$$
Cela veut dire, que
$$
w'\thicksim -\frac{1}{2}\{ y,u\} \thicksim -\frac{1}{2}\{y, u-y^2\}
\thicksim -\frac{1}{2}\{ y, u'\}
\eqno(2.9.21)$$
dans l'homologie de $DHABD^{odd}({\Bbb Q})$ (o\`u $y=[x_1x_2]$;
$u=[x_1x_3]\cdot [x_2x_4]$;
$u'=[x_1x_3]\cdot [x_2x_4] -[x_1x_2]\cdot [x_3x_4]$). $w'$ peut 
donc \^etre obtenu des diagrammes de cordes. $\Box$

Dans le bidegr\'e $(3,6)$ on a les relations suivantes entre les 
$c_i,$
$i=1,\dots,5$:
$$
\left.\begin{array}{l}
c_1=0;\\
c_2=c_3=c_5;\\
c_4=2c_2.
\end{array}\right\}
\eqno(2.9.22)$$

Cela donne encore un g\'en\'erateur pair (primitif).

\vspace{4mm}

\noindent{\bf 3. L'homologie des complexes de Hochschild
pour les op\'erades et l'homologie des espaces de n\oe uds}

\vspace{2mm}

\noindent{\bf 3.1. Structure de superalg\`ebre de Lie sur une op\'erade 
lin\'eaire}

\vspace{2mm}

Soit $\{ {\cal O}(n), n\ge 0\}$ une op\'erade dans la cat\'egorie des espaces
vectoriels ${\Bbb Z}_2$-gradu\'es (ou $\Bbb Z$-gradu\'es). Cet objet est 
li\'e \`a l'ensemble des espaces $\{{\rm Hom}(V^{\otimes n},V),n\ge 0\}$
de toutes les op\'erations  $n$-aires sur l'espace vectoriel $V$, de la m\^eme
fa\c{c}on qu'un groupe alg\'ebrique est li\'e 
au groupe $GL(V)$ des transformations lin\'eaires. Une structure 
d'op\'erade sur la famille des espaces ${\Bbb Z}_2$-(ou $\Bbb Z$-)gradu\'es 
${\cal O}(n)$ consiste en une famille des op\'erations conservant la 
graduation
$$
\gamma :{\cal O}(k)
\otimes\bigl(
{\cal O}(n_1)\otimes\dots\otimes{\cal O}(n_k)\bigr)
\to{\cal O}(n_1+...+n_k)
\eqno(3.1.1)$$
qui satisfont certains axiomes naturels de l'associativit\'e.
Ces axiomes peuvent \^etre obtenus du mod\`ele de l'op\'erade
des homomorphismes  $\{{\rm Hom}(V^{\otimes n},V),n\ge 0\}$; 
ici $\gamma$ est la 
substitution des valeurs des $k$ op\'erations dans une op\'eration $k$-aire.
La d\'efinition pr\'ecise peut \^etre trouv\'ee dans 
[GJ], [GK], ou [KSV]. On sous-entend \'egalement une action du groupe 
sym\'etrique
$S_n$ sur chaque ${\cal O}(n)$ et un \'el\'ement de l'unit\'e 
$ id\in {\cal O}(1)$, analogues \`a l'action de $S_n$ sur 
${\rm Hom}(V^{\otimes n},V)$  par les transpositions des entr\'ees et \`a 
l'op\'erateur identit\'e dans ${\rm Hom}(V,V)$, respectivement.

Consid\'erons l'espace vectoriel ${\cal O}=\oplus_n{\cal O}(n)$ --- 
la somme de toutes les composantes de l'op\'erade. On va d\'esigner
par $\tilde x$ la parit\'e (ou le degr\'e) d'un \'el\'ement $x\in
{\cal O}(n)$. D\'esignons \'egalement par
$$
|x|:=\tilde x+n-1,
\eqno(3.1.2)$$
o\`u ``$n$'' correspond \`a $n$ entr\'ees, ``-1'' \`a une sortie
de l'\'el\'ement $x$ de l'op\'erade. Notons que les applications (3.1.1) 
conservent la graduation $|\,.\,|$.

D\'efinissons une collection d'op\'erations polylin\'eaires sur ${\cal O}$:
$$
x \{x_1,\dots,x_n\}:=\sum (-1)^{\epsilon}\gamma(x;id,\dots,
id,x_1,id,\dots,id,x_n,id,\dots,id )
\eqno(3.1.3)
$$
pour $x,x_1,\dots,x_n\in{\cal O}$, o\`u on prend la somme par toutes
les substitutions possibles de $x_1,\dots,x_n$ dans $x$ dans l'ordre
prescrit; $\epsilon:=\sum_{p=1}^n|x_p|i_p$, $i_p$ \'etant le nombre des 
entr\'ees avant $x_p$. Les op\'erations
$x\{x_1,\dots,x_n\}$ conservent (la graduation initiale et) la graduation
$|\,.\,|$, {\em i.e.}
$$|x\{x_1,\dots,x_n\}|=|x|+|x_1|+\dots +|x_n|.
\eqno(3.1.4)$$

D\'efinissons aussi
$$
x\{\}:=x,\eqno(3.1.5)$$
$$
x\circ y:=x\{y\}.\eqno(3.1.6)$$

On peut facilement v\'erifier l'identit\'e suivante
$$   
x \{x_1,\dots,x_n\}\{y_1,\dots,y_n\}=$$
$$\sum_{0\le i_1\le j_1\le\dots\le i_m\le j_m\le n}
(-1)^{\epsilon}x\{y_1,\dots,y_{i_1},x_1\{y_{i_1+1},\dots,y_{j_1}\},
y_{j_1+1},\dots,y_{i_m},$$
$$
x_m\{y_{i_m+1},\dots,y_{j_m}\},y_{j_m+1},\dots,
y_n\},
\eqno(3.1.7)$$
o\`u $\epsilon:=\sum_{p=1}^m\left(|x_p|\sum_{q=1}^{i_p}|y_q|\right),$ 
autrement dit
le signe est obtenu, quand on fait passer les $x_i$ par les $y_j$.

{\bf Remarque 3.1.8.} L'identit\'e (3.1.7) pour le cas $m=n=1$ entra\^\i ne,
que le crochet
$$
[x,y]:=x\circ y-(-1)^{|x||y|}y\circ x
\eqno(3.1.9)$$
d\'efinit une structure de 
superalg\`ebre de Lie ${\Bbb Z}_2$-(ou $\Bbb Z$-)gradu\'ee
sur ${\cal O}$. $\Box$

\vspace{4mm}

\noindent{\bf 3.2. Complexes de Hochschild}

\vspace{2mm}

Soient ${\cal O}=\oplus_{n=0}^{+\infty} {\cal O}(n)$ une op\'erade 
quelconque,
${\cal AS}=\oplus_{n=0}^{+\infty} {\cal AS}(n)$ l'op\'erade
des alg\`ebres associatives, alors tout morphisme (dans la cat\'egorie
des op\'erades)
$$
\Pi:{\cal AS}\longrightarrow {\cal O}
\eqno(3.2.1)$$
d\'efinit un \'el\'ement $\Pi (m_2)\in{\cal O}(2)$, o\`u $m_2\in {\cal AS}(2)$
correspond \`a l'op\'eration de multiplication.

Notons que $m_2\circ m_2=0$ (c'est la condition de l'associativit\'e de
la multiplication). On a donc
$
[\Pi(m_2),\Pi(m_2)]=2\Pi(m_2\circ m_2)=0.$
Alors ${\cal O}$ devient une superalg\`ebre de Lie
${\Bbb Z}_2$-(ou $\Bbb Z$-)gradu\'ee
diff\'erentielle avec la diff\'erentielle
$$
\partial x=[\Pi (m_2),x]=\Pi (m_2)\circ x-(-1)^{|x|}x\circ\Pi (m_2).
\eqno(3.2.2)$$

Le complexe $({\cal O},\partial )$ ainsi obtenu sera appel\'e {\it 
complexe de Hochschild} pour l'op\'erade ${\cal O}$ (munie d'un morphisme 
$\Pi: {\cal AS}\to {\cal O}$). Si ${\cal O}$ est une op\'erade des 
homomorphismes, alors ce complexe est le complexe de Hochschild (classique) 
pour une alg\`ebre associative.

On peut montrer, que ${\cal O}$ poss\`ede une structure d'alg\`ebre 
associative diff\'erentielle ${\Bbb Z}_2$-(ou $\Bbb Z$-)gradu\'ee 
avec la m\^eme diff\'erentielle $\partial$ et avec la multiplication
$$
x*y :=(-1)^{|x|+1}\Pi (m_2)\{x,y\}.
\eqno(3.2.3)$$

{\bf Th\'eor\`eme 3.2.4.} {\it La multiplication $*$ et le crochet $[.,.]$
induisent une structure d'alg\`ebre de Gerstenhaber
(voir la d\'efinition 2.5.3) dans l'homologie du complexe de Hochschild 
$({\cal O},\partial )$.} $\Box$

{\bf D\'emonstration du Th\'eor\`eme 3.2.4:} Le Th\'eor\`eme se d\'eduit
des formules suivantes d'homotopie:
$$
x *y -(-1)^{(|x|+1)(|y|+1)}y *x=(-1)^{|x|}(\partial (x\circ y)-
\partial x\circ y-(-1)^{|x|}x\circ\partial y).
\eqno(3.2.5)$$
Cette formule entra\^\i ne la supercommutativit\'e de la multiplication $*$
en homologie.

$$
[x,y *z]-[x,y]*z-(-1)^{|x|(|y|+1)}y*[x,z]=$$
$$
(-1)^{|x|+|y|+1}(\partial(x\{y,z\})-(\partial x)\{y,z\}
-(-1)^{|x|}x\{\partial y,z\}-(-1)^{|x|+|y|}x\{y,\partial z\}).
\eqno(3.2.6)$$
Cette formule entra\^\i ne la compatibilit\'e (en homologie) entre le  
crochet $[.,.]$ et la 
multiplication $*$. $\Box$

\newpage

\noindent{\bf 3.3. Op\'erades des alg\`ebres de Poisson, Gerstenhaber,
 Batalin-Vilkovissky}

\vspace{2mm}

\noindent{\bf 3.3.1. Op\'erades des petites boules et complexes 
des $B$-diagrammes g\'en\'eralis\'es}

\vspace{2mm}

De mani\`ere analogue (au cas lin\'eaire) on peut consid\'erer les op\'erades
topologiques, c'est-\`a-dire des collections 
$\{ {\cal O}(n), n\ge 0\}$ d'espaces topologiques, o\`u
  le groupe sym\'etrique $S_n$ agit sur chaque
${\cal O}(n)$, et on fixe les
applications de composition
$$
\gamma :{\cal O}(k)\times\bigl(
{\cal O}(n_1)\times...\times{\cal O}(n_k)\bigr)
\to{\cal O}(n_1+...+n_k),
\eqno(3.3.1)$$
et un \'el\'ement unitaire $ id\in {\cal O}(1)$. 
On sous-entend de fa\c{c}on analogue, que l'op\'erade topologique satisfait
certaines conditions de l'associativit\'e et de la concordance avec l'action 
du groupe sym\'etrique.

L'homologie (sur un corps) d'une op\'erade topologique forme une op\'erade
lin\'eaire $\Bbb Z$-gradu\'ee.

Comme  exemple d'op\'erade topologique on peut prendre l'op\'erade
${\cal B}_d=\{{\cal B}_d(n), n\ge 0\}$ des petites boules de dimension
$d\ge 1$, voir [Co]. Les espaces ${\cal B}_d(n)$ sont les espaces de toutes 
les configurations de $n$ boules ordonn\'ees et disjointes dans une boules 
unitaire.

Les op\'erations de composition (3.3.1) sont des ``substitutions''
de $k$ configurations respectivement de $n_1,\dots,n_k$ boules
dans une configuration de $k$ boules.

L'homologie de l'op\'erade ${\cal B}_d$ d\'efinit les op\'erades 
lin\'eaires $\Bbb Z$-gradu\'ees suivantes (voir [Co]):

1)Pour $d=1$ cela donne l'op\'erade des alg\`ebres associatives.

2)Pour $d=2k+1$, $k=1,2,\dots$, --- l'op\'erade des alg\`ebres de Poisson
avec un crochet de degr\'e $1-d$.

3)Pour $d=2k$, $k=1,2,\dots$, --- l'op\'erade des alg\`ebres de Gerstenhaber
avec un crochet de degr\'e $1-d$.

{\bf Remarque 3.3.2.} L'op\'erade des alg\`ebres associatives a une filtration
naturelle; le quotient gradu\'e correspondant est l'op\'erade des 
alg\`ebres de Poisson avec un crochet de degr\'e 0. $\Box$

{\bf Th\'eor\`eme 3.3.3.} {\it Le complexe de Hochschild pour 
l'op\'erade des alg\`ebres de Poisson (resp. de Gerstenhaber) est 
naturellement isomorphe au complexe des $B$-diagrammes g\'en\'eralis\'es
(voir la section 1.3) pour $d$ impair (resp. pair).} $\Box$

{\bf D\'emonstration du Th\'eor\`eme 3.3.3:} La seule difficult\'e
technique est
 de bien comprendre ce qui se passe avec les signes. $\Box$

Le lien entre l'op\'erade des petites boules et les espaces des n\oe uds
n'est pas fortuit. Il est facile de voir, que l'espace ${\cal B}_d(n)$
des configurations de $n$ boules ordonn\'ees et disjointes
dans la boule unitaire est homotopiquement \'equivalent \`a 
l'espace des applications injectives de l'ensemble fini
(=vari\'et\'e de dimension z\'ero) $\{ 1,\dots,n\}$ dans ${\Bbb R}^d$,
ce dernier espace-l\`a peut \^etre vu comme un analogue de 
dimension finie de l'espace
des n\oe uds.

\vspace{2mm}

\noindent{\bf 3.3.2. L'op\'erade des alg\`ebres de Batalin-Vilkovissky}

{\bf D\'efinition 3.3.4.} Une alg\`ebre de Gerstenhaber ${\cal A}$ est 
appel\'ee {\it alg\`ebre de Batalin-Vilkovissky}, si sur ${\cal A}$
agit une application impaire
$$
\delta :{\cal A}\longrightarrow {\cal A},
\eqno(3.3.5)$$
telle que

1)$\delta^2=0$,

2)$\delta (a\cdot b)=\delta(a)\cdot b +(-1)^{\tilde{a}}a\cdot \delta(b)
+(-1)^{\tilde{a}-1}[a,b]$. $\Box$

1) et 2) entra\^\i nent

3)$\delta([a,b])=[\delta(a),b]+(-1)^{\tilde{a}-1}[a,\delta(b)].$

Si l'on est dans la cat\'egorie des espaces $\Bbb Z$-gradu\'es, alors
on d\'emande, que $\delta$ soit du m\^eme degr\'e que le crochet
$[.,.]$.

{\bf Th\'eor\`eme 3.3.6.} {\it Le complexe de Hochschild pour l'op\'erade
des alg\`ebres de Batalin-Vilkovissky est naturellement
isomorphe au complexe des $B_*$-diagrammes g\'en\'eralis\'es
pour $d$ pair.} $\Box$

{\bf D\'emonstration du Th\'eor\`eme 3.3.6:} Analogue \`a celle du 
Th\'eor\`eme 3.3.3. $\Box$

L'op\'erade des alg\`ebres de Batalin-Vilkovissky a aussi une
interpr\'etation g\'eom\'etrique comme l'homologie d'une certaine op\'erade 
topologique.

Consid\'erons l'op\'erade des petits disques (boules de dimension $d=2$).
Avec une confi\-guration on va \'egalement fixer un point sur le bord de chacun
des disques de la configuration. L'espace ${\cal B}'_d(n)$ ainsi obtenu
des configurations des disques ``marqu\'ees'' est \'evidemment
${\cal B}_d(n)\times\underbrace{S^1\times\dots\times S^1}_{n}$.
Les op\'erations de composition (3.3.1) sont d\'efinies comme substitution 
de $k$ configurations respectivement de $n_1,\dots,n_k$ disques ``marqu\'es''
dans une configuration de $k$ disques ``marqu\'es'';
sous cette substitution on tourne chacune des $k$ configurations par l'angle
correpondant au point (sur le bord) marqu\'e, en supposant que le point 
fix\'e du disque unitaire est toujours, par exemple, le point (1,0).

\begin{center}
{\bf R\'ef\'erences}
\end{center}
\begin{itemize}

\item[{[}BN{]}]D.Bar-Natan, {\it On the Vassiliev knot invariants},
Topology, 34 (1995), p.p. 423--472.

\item[{[}ChD{]}]S.V.Chmutov, S.V.Duzhin, {\it An upper bound for the number
of Vassiliev knot invariants}. J. of Knot Theory and its Ramifications 3
(1994), p.p.141-151.

\item[{[}ChDL{]}]S.V.Chmutov, S.V.Duzhin, S.K.Lando, {\it Vassiliev knot
invariants.}  I. Introduction, In: Singularities and Bifurcations,
Providence, RI: AMS, 1994, p.p. 117--126 (Adv. in Sov Math. 21).

\item[{[}Co{]}]F.R.Cohen, {\it The homology of $C_{n+1}$-spaces, $n\ge 0$,
The homology of iterated loop spaces}, Lecture Notes in Mathematics 533,
Springer--Verlag, Berlin, 1976, p.p.207--351.

\item[{[}G{]}]V.Ginzburg, {\it Resolution of diagonals 
and moduli spaces}, hep-th/9502013.

\item[{[}GJ{]}]E.Getzler, J.D.S.Jones, {\it Operads, homotopy algebra
and iterated integrals for double loop spaces}, preprint hep-th/9403055,
Departement of Mathematics, Massachusetts Institute of Technology, March 1994.

\item[{[}GK{]}]V.Ginzburg, M.Kapranov, {\it Koszul duality for operads},
Duke Math. J. 76 (1994), p.p.203--272.

\item[{[}GV{]}]M.Gerstenhaber, A.Voronov, {\it Homotopy $G$-algebras and 
moduli space operad}, Intern. Math. Res. Notices (1995), No.3, p.p. 141--153.
 
\item[{[}K1{]}]M.Kontsevich, {\it Vassiliev's knot invariants}, Adv. in Sov.
Math., vol.16, part 2, AMS, Providence, RI, 1993, p.p.137--150.

\item[{[}K2{]}]M.Kontsevich, {\it Private communication with V.Vassiliev},
April 1994, Texel Island.

\item[{[}Kn{]}]J.A.Kneissler, {\it The number of primitive Vassiliev
invariants up to degree twelve}, q-alg/97060222.

\item[{[}KSV{]}]T.Kimura, J.Stasheff, A.A.Voronov, {\it On operad structures
of moduli spaces and string theory}, to appear in Comm. Math. Phys.

\item[{[}L{]}]S.K.Lando. {\it On primitive elements in the bialgebra of chord
diagrams}, AMS Translations (2), vol.180, 1997, p.p.167--174.

\item[{[}NS{]}]K.Y.Ny, T.Stanford, {\it On Gusarov's groups of knot.} To
appear in Proc. Camb. Phil. Soc.

\item[{[}P{]}]F.Patras, {\it L'alg\`ebre
des descentes d'une big\`ebre gradu\'ee}, J. of Algebra. 170, 2 (1994), 
p.p.547--566.

\item[{[}S{]}]A.Stoimenow, {\it Enumeration of chord diagrams and an upper
bound for Vassiliev inva\-riants.} J. of Knot Theory and its Ramifications 7
(1998) 93-114.

\item[{[}T{]}]V.Tourtchine, {\it Sur les questions combinatoires de la 
th\'eorie spectrale des n\oe uds}, PHD thesis, Universit\'e Paris 7, 
(June 2001) to appear.

\item[{[}V1{]}]V.A.Vassiliev, {\it Cohomology of knot spaces.} In: Adv. in
Sov.  Math.; Theory of Singularities and its Applications (ed. V.I.Arnol'd).
AMS, Providence, R.I., 1990, p.p.23--69.

\item[{[}V2{]}]V.A.Vassiliev. {\it Stable homotopy type of the complement to
affine plane arrangements.} Preprint 1991.

\item[{[}V3{]}]V.A.Vassiliev. {\it Complexes of connected graphs}. In:
Gelfand's Mathematical Seminars, 1990--1992. L.Corwin, I.Gelfand, J.Lepovsky,
eds.  Basel: Birkh\"{a}user, 1993, p.p.223--235.

\item[{[}V4{]}]V.A.Vassiliev. {\it Complements of Discriminants of Smooth
Maps:  Topology and Applications.} Revised ed. Providence, R.I.: AMS, 1994
(Translation of Mathem. Monographs, 98).

\item[{[}V5{]}]V.A.Vassiliev, {\it Topology of two-connected graphs and
homology of spaces of knots}, Preprint 1999.

\item[{[}Z{]}]D.Zagier, {\it Vassiliev invariants and a strange identity
related to the Dedekind eta-function}, preprint 1999.

\end{itemize}

\vspace{2mm}

\hspace{3cm} Victor Tourtchine

\hspace{3cm} Universit\'e Ind\'ependante de Moscou,

\hspace{3cm} Universit\'e de Paris 7

\hspace{3cm} Russia, 121002 Moscow,

\hspace{3cm} B.Vlassjevskij 11, MCCME 

\hspace{3cm} e-mail: turchin@mccme.ru, tourtchi@acacia.ens.fr

\end{document}